\documentclass[10pt]{amsart}
\usepackage{mathptmx}
\usepackage{amsmath}
\usepackage{amssymb}
\usepackage{array}
\usepackage{geometry}
\usepackage[bookmarks=true,colorlinks=true, pdfstartview=FitV, linkcolor=black, citecolor=blue, urlcolor=black]{hyperref}

\usepackage{color}
\definecolor{DarkRed}{rgb}{0.55,.00,0.2}
\definecolor{DarkGrey}{rgb}{0.35,.35,0.35}

\theoremstyle{definition}

\theoremstyle{remark}

\numberwithin{equation}{section}

\begin{document}

\title{ On the generalized Dirichlet beta and Riemann zeta functions and Ramanujan-type formulae for beta and zeta values }

\author{S. Yakubovich}
\address{Department of Mathematics, Faculty of Sciences,  University of Porto,  Campo Alegre str.,  687; 4169-007 Porto,  Portugal}
\email{ syakubov@fc.up.pt}


\date{\today}
\maketitle

\markboth{\rm \centerline{ S.  Yakubovich}}{}

\markright{\rm \centerline{Generalized Dirichlet beta function}}

\begin{abstract}  We define the generalized Dirichlet beta  and Riemann zeta functions in terms of the   integrals, involving powers of the hyperbolic secant and cosecant functions. The corresponding functional equations are  established.   Some consequences of the Ramanujan identity for zeta values at odd integers are investigated and new formulae of the Ramanujan type are obtained.  \end{abstract}

\section{Introduction and preliminary results}

Throughout the text, $\mathbb{N}$ will denote the set of all positive integers, $\mathbb{N}_0 = \mathbb{N} \cup \{0\}$, whereas $\mathbb{R}$ and $\mathbb{C}$ the field of the real and complex numbers, respectively. The subject of the present investigation is the following elementary integrals of two complex variables $(z,s) \in \mathbb{C}\times\mathbb{C}$

$$\beta(z,s) = {1\over \Gamma(s)} \int_0^\infty {x^{s-1}\over \cosh^z x} dx,\quad {\rm Re}\ z > 0,\quad   {\rm Re}\ s > 0,\eqno(1.1)$$

$$\zeta(z,s) = {1\over \Gamma(s)} \int_0^\infty {x^{s-1}\over \sinh^z x} dx,\quad   {\rm Re}\  s >  {\rm Re}\ z > 0,\eqno(1.2)$$
where $\Gamma(s)$ is Euler's gamma function [14].  It is easily seen that integral (1.1) converges absolutely and uniformly by each variable on the sets $ {\rm Re}\ z \ge z_0 >0,\  {\rm Re}\ s \ge s_0 >0$. Therefore $\beta(z,s)$ represents an analytic function by each variable $z, s$ in the region $D \equiv \{ (z,s) \in \mathbb{C}\times\mathbb{C}:\  {\rm Re}\ z > 0,\quad   {\rm Re}\ s > 0\}$. Consequently, by virtue of the classical Hartogs theorem it is analytic in $D$ as a function of two complex variables.  Similarly, it concerns the function $\zeta(z,s)$ which is analytic in the region $D_1 \equiv \{ (z,s) \in \mathbb{C}\times\mathbb{C}:\  {\rm Re}\ s >  {\rm Re}\ z > 0\}$.  We will call integrals (1.1), (1.2) as the generalized Dirichlet beta and Riemann zeta functions, respectively. 

When $z=1$, we have via Entry 2.4.3.2 in [9], Vol. I $\beta (1, s) = 2\beta(s),$ where $\beta (s)$ is the Dirichlet beta function

$$ \beta (s) = \sum_{n=0}^\infty {(-1)^n\over (2n+1)^s},\eqno(1.3)$$
involving the familiar Catalan constant  $G= \beta(2)$ and satisfying the following functional equation

$$ \beta(1-s) = \left({\pi\over 2}\right)^{-s} \Gamma(s)  \sin\left({\pi s\over 2}\right) \beta (s),\quad s \in \mathbb{C}.\eqno(1.4)$$
The case $z=2$ exhibits the Riemann zeta function $\zeta(s)$ (cf. Entry 2.4.3.3 in [6], Vol. I), namely $\beta(2,s) = 2^{2-s} \left( 1- 2^{2-s} \right)  \zeta (s-1),$ where

$$ \zeta (s) = \sum_{n=0}^\infty {1\over n^s},\quad {\rm Re}\ s > 1\eqno(1.5)$$
which obeys  the functional function 

$$\zeta(s)= 2^s \pi^{s-1} \Gamma(1-s)  \sin\left({\pi s\over 2}\right) \zeta (1-s),\quad s \in \mathbb{C}.\eqno(1.6)$$
Accordingly, the values  $\zeta(1,s) = \left(2- 2^{1-s}\right) \zeta(s),\  \zeta(2,s) = 2^{2-s} \zeta(s-1)$ are due to Entries 2.4.3.1. and 2.4.3.3 in  [6], Vol. I.

In the sequel we will employ the sequence of polynomials $\{p_n(x;\alpha)\}_{n\ge 0}$ [7] of exactly degree $n$ written in the operator form

$$p_n(x;\alpha) = (-1)^n e^x x^{-\alpha}\mathcal{A}^n e^{-x} x^\alpha,\quad n \in \mathbb{N}_0,\   {\rm Re}\ \alpha  > - {1\over 2},\eqno(1.7)$$
where $\mathcal{A}$ is the Bessel second order differential operator 

$$\mathcal{A} = x^2- x {d\over dx} x {d\over dx},\eqno(1.8)$$
whose eigenfunction is the modified Bessel function of the third kind $K_{i\tau}(x)$ of the argument $x >0$ and pure imaginary subscript $i\tau$, i.e.

$$\mathcal{A} K_{i\tau}(x) = \tau^2 K_{i\tau}(x).\eqno(1.9)$$
This function can be defined by the Fourier integral 

$$K_{i\tau}(x) = \int_0^\infty e^{-x\cosh u}\cos(\tau u) du,\quad x >0,\ \tau \in \mathbb{R}\eqno(1.10)$$
and it is the kernel of the Kontorovich-Lebedev transform [13]

$$(KL f)(\tau)= \int_0^\infty K_{i\tau}(x) f(x) dx\eqno(1.11)$$
for a suitable class of functions and with the reciprocal inverse transform

$$f(x)= {2\over \pi^2 x}   \int_0^\infty \tau \sinh(\pi\tau) K_{i\tau}(x) (KL f) (\tau) d\tau,\eqno(1.12)$$
where the convergence sense of the integrals (1.11), (1.12) is presumed accordingly.   Note, that the polynomial sequence  $\{p_n(x; 0)\}_{n\ge 0}$ is considered in [12].     The explicit expression for the sequence $\{p_n(x; \alpha)\}_{n\ge 0}$ is establish in [7], and we have the formula 

$$ p_n(x; \alpha) = \sum_{\nu =0}^n \bigg\{ {2^{\nu+1}\over \nu!} \left(\alpha + {1\over 2}\right)_\nu \ \sum_{\mu=0}^\nu \binom{\nu}{\mu} {(-1)^\mu 
\Gamma(2\alpha + \mu)\over \Gamma(2\alpha+\mu+\nu+1)} (\alpha+\mu)^{2n+1} \bigg\} x^\nu,\quad n \ge 0,\eqno(1.13)$$
where $ {\rm Re} \alpha > -1/2$ and

$$(a)_\nu= a(a+1)\dots (a+\nu-1)= {\Gamma(a+\nu)\over \Gamma(a)} \eqno(1.14)$$ 
is the Pochhammer symbol. Returning to the definition of the generalized Dirichlet beta function (1.1),  we can immediately expand it in terms of the Dirichlet series. Namely, it has under the additional condition ${\rm Re}\ (s -z ) > 0$  

$$ \beta(z,s) = {2^z \over \Gamma(s)} \int_0^\infty {x^{s-1} e^{-zx} \over (1+ e^{-2x})^z } dx = {2^z \over \Gamma(s)} \int_0^\infty x^{s-1}  e^{-zx} \sum_{n=0}^\infty {(-1)^n (z)_n \over n!} e^{-2nx} $$

$$= 2^z \sum_{n=0}^\infty {(-1)^n\  (z)_n \over n! \ (2n+z)^s},$$
i.e.

$$ \beta(z,s) =  2^z \sum_{n=0}^\infty {(-1)^n\  (z)_n \over n! \ (2n+z)^s},\quad {\rm Re}\ (s -z ) > 0,\eqno(1.15)$$
where the interchange of the order of integration and summation is allowed due to the absolute convergence  

$$\int_0^\infty x^{{\rm Re} s-1}  e^{-  x {\rm Re}\ z} \sum_{n=0}^\infty {({\rm Re}\ z)_n \over n!} e^{-2nx} < \infty, \quad {\rm Re}\ (s -z ) > 0.$$
In the same manner the series expansion is obtained for the generalized Riemann zeta function, and we have 

$$ \zeta(z,s) =  2^z \sum_{n=0}^\infty { (z)_n \over n! \ (2n+z)^s},\quad {\rm Re}\ (s -z ) > 0.\eqno(1.16)$$
Other two expansions  follow immediately from definitions (1.1), (1.2) and  the interchange of the order of integration and summation. Precisely, it has

$$ \beta (z,s) = 2^{z/2 -s}\  \sum_{k=0}^\infty (-1)^k \ {(z/2)_k \over k!} \ \beta \left({z\over 2} +k, s\right),$$

$$ \zeta (z,s) = 2^{z/2 -s}\  \sum_{k=0}^\infty \ {(z/2)_k \over k!} \ \beta \left({z\over 2} +k, s\right).$$

Further properties of these functions  and their  values at integers will be established with the use  of the central factorial numbers [10] of the first kind $t(n,\nu)$ and the second kind $T(n,\nu).$  They satisfy the following triangular relations 

$$t(n,\nu)= t(n-2,\nu-2) - {1\over 4} (n-2)^2 t(n-2,\nu),\ \quad   0\le \nu \le n,\eqno(1.17)$$

$$t(n,0)= t(0,n)= \delta_{n,0},\quad t(n,n)=1, \quad n \ge 0,\eqno(1.18)$$

$$T(n,\nu)= T(n-2,\nu-2) - {1\over 4} \nu^2 T(n-2,\nu),\ \quad   0\le \nu \le n,\eqno(1.19)$$

$$T(n,0)= T(0,n)= \delta_{n,0},\quad T(n,n)= 1, \quad n \ge 0,\eqno(1.20)$$
where $\delta_{n,\nu}$ is the Kronecker delta.   When $\nu \ge n+1$ or $(-1)^n + (-1)^\nu =0$, it has necessarily, $t(n,\nu) = T(n,\nu) = 0$.  Splitting our analysis into the cases of even or odd order, i.e.

$$t_E(n,\nu) := t(2n,2\nu), \quad  T_E(n,\nu) := T(2n,2\nu),$$

$$ t_O(n,\nu) := t(2n+1,2\nu+1), \quad  T_O(n,\nu) := T(2n+1,2\nu+1),\  0\le \nu \le n,\eqno(1.21)$$
we mention four reciprocities which will be employed in the sequel

$$ \prod_{k=1}^{n} \left[ k^2+ x^2\right] = \sum_{k=0}^n (-1)^{n+k} t_E(n+1, k+1) x^{2k},\eqno(1.22)$$

$$ x^{2n} = \sum_{k=0}^n (-1)^{n+k} T_E(n+1, k+1) \prod_{q=1}^{k} \left[ q^2+ x^2\right],\eqno(1.23)$$

$$ \prod_{k=1}^{n} \left[ \left(k- {1\over 2}\right)^2+ x^2\right] = \sum_{k=0}^n (-1)^{n+k} t_O(n, k) x^{2k},\eqno(1.24)$$

$$ x^{2n} = \sum_{k=0}^n (-1)^{n+k} T_O(n, k) \prod_{q=1}^{k} \left[ \left(q- {1\over 2} \right)^2+ x^2\right].\eqno(1.25)$$
Closed-form expression for the central factorial numbers of the second kind is given by the formula (cf. [10])

$$ T(n,\nu)= {1\over \nu!} \sum_{\mu=0}^\nu \binom{\nu}{\mu} (-1)^\mu \left({\nu\over 2} - \mu\right)^n,\quad   0\le \nu \le n, \quad n,\nu \in \mathbb{N}_0.\eqno(1.26)$$
Further, in the sequel we will appeal to the reciprocal formulae of the Mellin transform [14]

$$f^*(s)= \int_0^\infty f(x) x^{s-1} dx,\eqno(1.27)$$

$$f(x)= {1\over 2\pi i} \int_{\gamma - i\infty}^{\gamma + i\infty} f^*(s) x^{-s} ds,\eqno(1.28)$$
employing them under suitable conditions to generalize the remarkable Ramanujan formula for zeta values at odd integers (cf. [1])

$$a^{-k}\left( {1\over 2} \ \zeta (2k+1) + \sum_{n=1}^\infty {n^{-2k-1} \over e^{2a n} -1} \right) =  (- b)^{-k}\left( {1\over 2} \ \zeta (2k+1) + \sum_{n=1}^\infty {n^{-2k-1} \over e^{2 b n} -1} \right) $$

$$- 2^{2k} \sum_{j=0}^{k+1} {(-1)^j B_{2j} B_{2k+2-2j}\over (2j)!(2k+2-2j)!}  a^{k+1-j}  b^j,\quad  a, b > 0, \ ab = \pi^2,\ k \in \mathbb{N}.\eqno(1.29)$$
As far as we aware,  the method of the Mellin transform and Mellin-Barnes type integrals  originated from [8], [2].  We note also seminal papers [3], [4] and some rapidly convergent series for the Riemann zeta-values being obtained via the modular relation in [5].   We will apply various Mellin-Barnes type integral representations for the hyperbolic functions, manipulating with functional equations (1.4), (1.6) for the Dirichlet beta and Riemann zeta functions and differentiation under the integral sign with respect to parameters.   In particular, it will be proved the following Ramanujan-type formula 

$$  (-1)^{k}  \binom{2k-1+m}{m}  \left({\pi\over a}\right)^{2k}\  \zeta(2k+1) $$

$$ +   (2\pi)^{2k+1} \sum_{q=0}^{k+1} (-1)^{k+q}  {B_{2q} B_{2(k+1-q)} \over  (2q)!  (2(k+1-q))!} \  \binom{2(q-1)+m}{m}  \left({\pi\over a}\right)^{2q-1}$$

$$= 2 {b^{m}\over m!}   {d^{m}\over d b^{m}}  \sum_{n=1}^\infty {n^{-2k-1} \over  e^{2 b n} -1 } +    2 \sum_{r=1}^{m-1}  {r\over (r-1)!} \ \binom{m-1}{r}  \ b^{r} \  {d^{r}\over d b^{r}}  \sum_{n=1}^\infty {n^{-2k-1} \over  e^{2 b n} -1 } $$

$$ + 2  (-1)^k \  \left({\pi\over a}\right)^{2k}    \sum_{r=0}^{m} {(-1)^{r+1}\over r!}  \binom{2k-1+m-r}{2k-1}  a^{r}  {d^r\over d a^r}  \sum_{n=1}^\infty {n^{-2k-1} \over  e^{2 a n} -1 },\eqno(1.30)$$
where $a, b > 0$ with $ab = \pi^2,\  k, m \in \mathbb{N}$ and $B_n$ are Bernoulli numbers (cf. [9], Vol. I).  Denoting by $F_k(x)$ the series 

$$F_{k}(x) =  \sum_{n=1}^\infty {n^{- k} \over  e^{2  n x} -1 },\quad x > 0,\ k \in \mathbb{N},\eqno(1.31)$$
 and by $F_{k}^{(j)}(x)$ the corresponding derivative of the order $j$, we rewrite formula (1.30) for the case $a=b=\pi$ to obtain for $k, m \in \mathbb{N}$

$$   \binom{2k-1+m}{m} \bigg[ \zeta(2k+1) + 2 \sum_{n=1}^\infty {n^{-2k-1} \over  e^{2  n x} -1 } \bigg] +   (2\pi)^{2k+1} \sum_{q=0}^{k+1}   { (-1)^{q} \ B_{2q} B_{2(k+1-q)} \over  (2q)!  (2(k+1-q))!} \  \binom{2(q-1)+m}{m}$$

$$=  2 \left( (-1)^k - (-1)^m \right)  {\pi^{m}\over m!} F_{2k+1}^{(m)}(\pi)   +   2  \sum_{r=1}^{m-1} \bigg[  (-1)^k  r^2 \ \binom{m-1}{r} + (-1)^{r+1} \binom{2k-1+m-r}{m-r} \bigg]  \ {\pi^{r}\over r!}  \  F_{2k+1}^{(r)}(\pi).\eqno(1.32)$$
Letting $m=1$, we obtain a companion for the Lerch formula [1] for even $k=2n,\ n \in \mathbb{N}$. Precisely, we find the identity

$$     \zeta(4n+1) + 2  \sum_{k=1}^\infty {k^{-4n-1} \over  e^{2 \pi k} -1 } + {\pi\over 2n} \sum_{k=1}^\infty {k^{-4n} \over  \sinh^2 (\pi k)} = {2^{4n}\over 2 n}\  \pi^{4n+1} \sum_{q=1}^{2n+1}  { (-1)^{q+1} \  B_{2q} B_{4n+2-2q}  \over  (2q-1)!  (4n+2-2q)!},\eqno(1.33)$$
where the equality 
$$  {2^{4n}\over 2 n}\  \pi^{4n+1} \sum_{q=0}^{2n+1}  { (-1)^{q+1} \ (2q-1)\  B_{2q} B_{4n+2-2q}  \over  (2q)!  (4n+2-2q)!} =  {2^{4n}\over 2 n}\  \pi^{4n+1} \sum_{q=1}^{2n+1}  { (-1)^{q+1} \  B_{2q} B_{4n+2-2q}  \over  (2q-1)!  (4n+2-2q)!}$$
is used  as a consequence of  the Ramanujan formula (1.29) for this case. In particular, $n=1$ yields the following result

$$     \zeta(5) + 2  \sum_{k=1}^\infty {1\over  k^5 (e^{2 \pi k} -1) } + {\pi\over 2} \sum_{k=1}^\infty {1\over  k^4 \sinh^2 (\pi k)} = {13\over 3780}\ \pi^5.\eqno(1.34)$$
Note, that the odd case $m=1,\ k=2n+1$  confirms the familiar Lerch formula [1] written in the form

$$    \zeta(4n+3) =  {2^{4n+2} \pi^{4n+3}\over 2n+1}  \sum_{q=0}^{2n+2}   { (-1)^{q+1} \   (2q-1)\ B_{2q} B_{4n+4- 2q} \over  (2q)!  (4n+4-2q)!} - 2 \sum_{k=1}^\infty {k^{- 4n-3} \over  e^{2 \pi k} -1 },\eqno(1.35)$$
and its very particular case

$$ \zeta(3) + 2  \sum_{k=1}^\infty {1\over  k^3 (e^{2 \pi k} -1) } = {7\over 180}\ \pi^3.\eqno(1.36)$$
In fact, this is because  it can be  easily seen the equality

$${1\over 2n+1}  \sum_{q=0}^{2n+2}   { (-1)^{q+1} \   (2q-1)\ B_{2q} B_{4n+4- 2q} \over  (2q)!  (4n+4-2q)!} =  \sum_{q=0}^{2n+2}   { (-1)^{q+1} \   B_{2q} B_{4n+4- 2q} \over  (2q)!  (4n+4-2q)!}.$$

\section{Properties of  $\beta(z,s)$}

We begin, establishing straightforward functional recurrence relation for integral (1.1).   Indeed,  assuming the conditions ${\rm Re}\ z,\  {\rm Re}\ s > 2$,  from the definition (1.1) of the generalized Dirichlet beta function and integration by parts twice we derive 

$$\beta (z,s) =  \beta (z-2, s) - {1\over \Gamma(s)} \int_0^\infty {x^{s-1} \sinh^2 x \over \cosh^z x} dx =   {z-2\over z-1}\  \beta (z-2, s) - {  \beta (z-2, s-2) \over (z-1)(z-2)},$$
where integrated terms are vanished. Hence, changing $z-2,\ s-2$ by $z, s$, we find the following recurrence relation 

$$\beta (z,s) =  \ z^2 \beta (z, s+2) - z(z+1) \beta (z+2, s+2),\eqno(2.1)$$
which can be analytically continued to the region $\{ (z,s) \in \mathbb{C}\times\mathbb{C}:\  {\rm Re}\ z, \  {\rm Re}\ s > -2\}.$  Incidentally, relation (2.1) can be generalized, employing the polynomial sequence (1.13) $\{p_n(x; \alpha)\}_{n\ge 0}$.  Actually, we have

{\bf Theorem 1}. { \it  Let $n \in \mathbb{N}_0$ and $c_{n,k}(\alpha)$ be polynomial coefficients of the sequence $\{p_n(x; \alpha)\}_{n\ge 0}$. Then the following equality takes place}

$$ \sum_{k=0}^n  { (z/2)_k \over   2^{k} }\  c_{n,k}\left({z\over 2} \right)\  \beta (z+2k, s+2n) =    4^{-n }   \beta (z, s),\  {\rm Re}\  z > 0,\   {\rm Re}\ s > -2n,\eqno(2.2)$$

\begin{proof}   Taking for now $0 < {\rm Re}\ s < 1$, we use the Parseval equality for the Fourier cosine transform [14] and Entries 2.5.3.10, 2.5.46.6 in [9], Vol. I to derive the representation

$$ \beta(z,s) = {2^{z-s}\over \pi \ \Gamma(z) }\ \cos\left({\pi s\over 2}\right) \int_0^\infty x^{-s}  \Gamma\left( {z\over 2}+ix\right)  \Gamma\left( {z\over 2} - ix\right) dx.\eqno(2.3)$$
When $z=1$, we know  that  $\beta (1, s) = 2\beta(s)$, and therefore

$$2\beta(s) = \cos\left({\pi s\over 2}\right) \int_0^\infty {x^{-s} \over \cosh (\pi x /2)}\ dx =  \left({2\over \pi} \right)^{1-s} \cos\left({\pi s\over 2}\right) \int_0^\infty {x^{-s} \over \cosh x}\ dx, $$
i.e.

$$ \beta(1-s) = \left({\pi\over 2}\right)^{-s} \Gamma(s)  \sin\left({\pi s\over 2}\right) \beta (s)$$
which is again  functional equation (1.4) for the Dirichlet beta function.  Meanwhile, Entry 2.16.6.4 in [9], Vol.II and the duplication formula for the gamma function suggest the representation from (2.3)

$$ \beta (z, s) = {2^{z/2 -s+1}\over \pi  \Gamma(z/2)}\ \cos\left({\pi s\over 2}\right) \  \int_0^\infty x^{-s}   \int_0^\infty t^{z/2-1} e^{-t} K_{ix}(t) dt dx.\eqno(2.4)$$
Hence, putting   $z+2k,\   k \in \mathbb{N}_0$ instead of $z$ and employing (1.13), equality (2.4) becomes

$$  \Gamma(z/2) \sum_{k=0}^n  c_{n,k}\left({z\over 2}\right) \  {(z/2)_k \over   2^{z/2 +k+1}}\    \beta (z+2k, s) =  {2^{ -s}\over \pi}\ \cos\left({\pi s\over 2}\right)  \int_0^\infty x^{-s}   \int_0^\infty t^{z/2-1} e^{-t}  p_n\left(t; {z\over 2}\right) K_{ix}(t) dt dx,\eqno(2.5)$$
where (see (1.13))

$$ c_{n,k}(\alpha) = {2^{k+1}\over k!} \left(\alpha + {1\over 2}\right)_k \ \sum_{\mu=0}^k \binom{k}{\mu} {(-1)^\mu 
\Gamma(2\alpha + \mu)\over \Gamma(2\alpha+\mu+k+1)} (\alpha+\mu)^{2n+1}.\eqno(2.6)$$
Now, employing (1.7), (1.9)  and integration by parts, we obtain 

$$   \Gamma(z/2) \sum_{k=0}^n  c_{n,k}\left({z\over 2}\right) \  {(z/2)_k \over   2^{z/2 +k+1}}\    \beta (z+2k, s) =  {2^{ -s} (-1)^n \over \pi}\ \cos\left({\pi s\over 2}\right)  \int_0^\infty x^{2n -s}   \int_0^\infty t^{z/2-1} e^{-t}  K_{ix}(t) dt dx.\eqno(2.7)$$
Thus,  recalling again Entry  2.16.6.4 in [9], Vol.II and comparing with (2.3), we observe 

$$ \sum_{k=0}^n  c_{n,k}\left({z\over 2}\right) \  {(z/2)_k \over   2^{k}}\    \beta (z+2k, s) =    {2^{1 -s } (-1)^n \over \sqrt \pi\  \Gamma((z+1)/2) }\ \cos\left({\pi s\over 2}\right)   \int_0^\infty x^{2n -s}\   \Gamma\left( {z\over 2} +ix\right)  \Gamma\left( {z\over 2}  - ix\right) dx$$

$$ = 4^{ -n}  \beta (z, s-2n).$$
Consequently, with a simple change of variable, we arrive at the relation (2.2) which can be, in turn,  analytically continued to the domain  ${\rm Re}\ s > - 2n$.  

\end{proof}

Further, let us consider the values of the generalized Dirichlet beta function (1.1) for positive even integers $z$.  Then, returning to (2.3),  we write

$$ \beta (2(n+1), s) = {2^{2n-s+2}\over \pi (2n+1)! }\ \cos\left({\pi s\over 2}\right) \int_0^\infty x^{-s}  \Gamma\left( n+1+ix\right)  \Gamma\left( n+1 - ix\right) dx$$

$$ = {2^{2n-s+2}\over  (2n+1)! }\ \cos\left({\pi s\over 2}\right) \int_0^\infty {x^{1-s}\over \sinh(\pi x)} \   \prod_{k=1}^{n} \left[ k^2+ x^2\right] dx,\quad  {\rm Re}\ s  <1.\eqno(2.8)$$
Meanwhile, the latter product can be expressed  in terms of the central factorial numbers of the first kind (see (1.22)). Hence,  appealing to Entry 2.4.3.1 in [9], Vol. I, we arrive at the equality

$$  \beta (2(n+1), s) = {2^{2n-s+3} (-1)^{n+1} \over  (2n+1)! }\ \cos\left({\pi s\over 2}\right) \sum_{k=1}^{n+1}  {(-1)^{k}\over \pi^{2k-s}}\   t_E(n+1, k) $$

$$\times \left( 1 - 2^{s-2k}\right) \Gamma(2k-s) \zeta (2k-s),\quad n \in \mathbb{N}_0$$ 
 which  can be  extended analytically to the whole complex plane, involving the functional equation (1.6) for the Riemann zeta function.   Thus it implies  the equality

$$  \beta (2n, s) = {4^{n-s} (-1)^{n} \over   (2n-1)! }\  \sum_{k=1}^{n}  t(2n, 2k) \left( 2^{2k} - 2^{s}\right) \zeta (1+s-2k),\   s \in \mathbb{C},\ n \in \mathbb{N}.\eqno(2.9)$$ 
In particular, if $s$ is an even integer $s= 2m,\ m \in \mathbb{Z}$, the generalized Dirichlet beta function (1.1)  has the values

$$  \beta(2n, 2m) = {4^{n-2m} (-1)^{n} \over   (2n-1)! }\  \sum_{k=1}^{n}  t(2n, 2k) \left( 4^{k} - 4^{m}\right) \zeta (1+2(m-k)),\quad m \in \mathbb{Z},\ n \in \mathbb{N}.\eqno(2.10)$$ 
Moreover, if $m$ is a nonpositive  integer, i.e. $m= -r,\ r \in \mathbb{N}_0$, we take into account  the corresponding values of the Riemann zeta function in terms of Bernoulli  numbers

$$\zeta (- n) = (-1)^n {B_{n+1}\over n+1},\quad n \in \mathbb{N}_0\eqno(2.11)$$
to derive the identity

$$  \beta(2n, - 2r) = {4^{2r-1} (-1)^{n+1} \over   (2n-1)! }\  \sum_{k=1}^{n}  t(2n, 2k) \left( 4^{k} - 4^{-r}\right) {B_{2(r+k)}\over r+k},\quad r, n \in \mathbb{N}.$$ 
It means that $ \beta(2(n+1), - 2r) \in \mathbb{Q}$.  When $m=1,\  n =0$ in (2.10), we see that the simple pole of the zeta function is compensated by a simple zero, and since $t(2,2)=1$, we find the value $\beta(2, 2) = \log 2$ which is confirmed via Entry 2.4.3.4 in [9], Vol. I.  Finally, for the general case $m \in \mathbb{N}$ we write

$$   \beta(2n, 2m) = {4^{n-m} (-1)^{n+1}  \over   (2n-1)! }\  \bigg[ \sum_{k=1}^{m-1}  t(2n, 2(m-k)) \left( 1- 2^{-2k} \right) \zeta (2k+1)\bigg.$$

$$+  t(2n, 2m) \log 2 -  {1\over 2}   \sum_{k=1}^{n-m}  t(2n, 2(k+m)) \left( 4^{k} - 1\right) {B_{2k}\over k}\bigg]  ,\quad m, n  \in \mathbb{N}.\eqno(2.12)$$ 
Taking $m=n$, we get, in particular,

$$   \beta(2n, 2n) = { (-1)^{n+1}  \over   (2n-1)! }\  \bigg[ \sum_{k=1}^{n-1}  t(2n, 2k) \left( 1- 2^{-2(n-k)} \right) \zeta (2(n-k)+1)+ \log 2 \bigg],\quad n \in \mathbb{N}.\eqno(2.13)$$
But if we let $m=n+2$, we use properties of the central factorial numbers [10] to write 

$$   \beta(2n, 2(n+2)) = { (-1)^{n+1}  \over   (2n-1)! }\   \sum_{k=0}^{n-1}  t(2n, 2(n-k)) \left( 4^{k+2} -  1 \right) {\zeta (2k+5)\over 4^{k+4}},\quad  n \in \mathbb{N}.\eqno(2.14)$$
Moreover, recalling the definition of the generalized Dirichlet beta function (1.1), it gives the values of the integrals  

$$   \int_0^\infty {x^{2n+3}\over \cosh^{2n} x } \ dx  =  (-1)^{n+1}  n(n+1)(2n+1) (2n+3)\   \sum_{k=0}^{n-1}  t(2n, 2(n-k)) \left( 4^{k+2} -  1 \right) {\zeta (2k+5)\over 4^{k+1}},\  n \in \mathbb{N}.\eqno(2.15)$$
In the meantime,  recalling  Entry 2.4.3.1 in [9], Vol. I and reciprocal identity (1.23), we have for $2n +1< {\rm Re}\ s < 2(n+1)$

$$\int_0^\infty {x^{s-1}\over \sinh(\pi x)}\ dx =  {2- 2^{1-s}\over \pi^s}\ \Gamma (s) \zeta (s) =  \sum_{k=0}^n (-1)^{n+k} T(2(n+1), 2(k+1)) \int_0^\infty {x^{s - 2n -1}\over \sinh(\pi x)}\  \prod_{q=1}^{k} \left[ q^2+ x^2\right]\ dx$$

$$=  {1 \over \pi} \sum_{k=0}^n (-1)^{n+k} T(2(n+1), 2(k+1)) \int_0^\infty x^{s - 2n -2} \  \Gamma\left( k+1+ ix \right)  \Gamma\left( k+1 - ix \right) dx,\ n \in \mathbb{N}_0.$$
Hence as above the Parseval equality for the Fourier cosine transform reduces to the equality

$$2^{s-2n} \ {1- 2^{-s}\over \pi^{s-1} }\ \Gamma (s) \zeta (s) =   \Gamma\left( s-2n-1\right) \sin\left({\pi s\over 2}\right) \sum_{k=0}^n {(-1)^{k} (2k+1)! \over 4^{k}}\  T(2(n+1), 2(k+1))  \int_0^\infty {x^{1+2n-s}\over  \cosh^{2(k+1)} x} \ dx.$$
Moreover, the functional equation (1.6) for the Riemann zeta function and the definition of the generalized Dirichlet beta function (1.1) allow to write it in the form

$$ 2^{s- 2n-2}  \left( 1-2^{s}\right) \ \zeta (1-s) =   \sum_{k=1}^{n+1} (-1)^{k+1}\ 4^{-k}   (2k-1)! \  T(2(n+1), 2k)\  \beta \left(2k, 2+2n-s\right),\eqno(2.16)$$
which can be continued analytically to all $s \in \mathbb{C}$.  In particular, letting $s=0$, we find

$$  \log 2  =   \sum_{k=1}^{n} (-1)^{k+1} 4^{n-k} (2k-1)!\  T(2n, 2k)\  \beta \left(2k, 2n\right),\quad n \in \mathbb{N}.\eqno(2.17)$$
Generally, for a nonnegative integer $s= m$ we write from (2.16)

$$ 2^{m- 2n-2}  \left( 1-2^{m}\right) \ \zeta (1-m) =   \sum_{k=1}^{n+1} {(-1)^{k+1}  (2k-1)!\over 4^k}\  T(2(n+1), 2k)\  \beta (2k,  2+2n-m),\eqno(2.18)$$
and, accordingly, for $s=2m$

$$4^{m- n-1}  \left( 1-4^{m}\right) \ \zeta (1-2m) =   \sum_{k=1}^{n+1} {(-1)^{k+1}  (2k-1)!\over 4^{k}}\  T(2(n+1), 2k)\  \beta (2k,  2(1+n-m)),\eqno(2.19)$$
but, in turn, for  $s=2m+1$ we have the formula
$$2^{2(m- n)-1}  \left( 1-2^{2m+1}\right) \ \zeta (-2m) =   \sum_{k=1}^{n+1} {(-1)^{k+1}  (2k-1)!\over 4^k}\  T(2(n+1), 2k)\  \beta (2 k,  1+2(n-m)).\eqno(2.20)$$
Hence, taking in mind (2.11), it gives

$$\sum_{k=1}^{n+1} (-1)^{k} 4^{-k} \  (2k-1)!\  T(2(n+1), 2k)\  \beta (2 k, 1+2(n-m)) = 0,\quad m \in \mathbb{N},$$

$$\sum_{k=1}^{n+1}  (-1)^{k+1} \ 4^{-k} \ (2k-1)! \  T(2(n+1), 2k)\  \beta (2 k, 2n+1) = 4^{-n-1},$$

$$\sum_{k=1}^{n+1}  (-1)^{k+1} 4^{-k} \  (2k-1)!\  T(2(n+1), 2k)\  \beta (2k,  2(1+n-m)) = 4^{m- n-1}  \left(4^{m}-1\right)  {B_{2m}\over 2m},\quad m \in \mathbb{N}.$$
However for nonpositive even and odd integers $s$ we have the following representations for the odd and even zeta values

$$  \left( 1-4^{-m}\right) \ \zeta (2m+1) =   \sum_{k=1}^{n}  (-1)^{k+1} 4^{m+n-k}  (2k-1)!\  T(2n, 2k)\  \beta (2k, 2(n+m)),\eqno(2.21)$$

$$   \left( 2- 4^{-m}\right) \ \zeta (2m) =   \sum_{k=0}^n (-1)^{k} 4^{1+m+n-k}  (2k+1)!\  T(2(n+1), 2(k+1))\  \beta (2(k+1),  3+2(n+m)).\eqno(2.22)$$
In particular, taking $m=1$, it gives

$$  \zeta (3) =   {16\over 3} \sum_{k=1}^{n}  (-1)^{k+1} 4^{n-k}  (2k-1)!\  T(2n, 2 k)\  \beta (2k, 2(n+1)),\quad n \in \mathbb{N}.$$

Returning to (2.16), we use the functional equation (1.6) in its left-hand side to write it in the form

$$ 2 \pi^{-s}  \left( 1-2^{s}\right) \ \cos\left({\pi s\over 2}\right) \Gamma(s) \zeta(s)  =   \sum_{k=1}^{n+1} (-1)^{k+1}\ 4^{n+1-k} \  (2k-1)! \  T(2(n+1), 2k)\  \beta (2k, \ 2+2n-s).$$
Hence subtracting two equalities, we find a kind of the functional equation for the generalized Dirichlet beta function

$$ \sum_{k=1}^{n} (-1)^{k+1} \ 4^{-k}  (2k-1)! \  T(2n, 2k)\bigg[   \pi^{-s}  \left( 1-2^{s}\right) \ \cos\left({\pi s\over 2}\right) \Gamma(s)   \beta (2k,  2n+s-1)\bigg.$$

$$\bigg. - 2^{-s}   \left( 1-2^{1-s}\right) \  \beta (2k,  2n-s)\bigg]= 0,\quad s \in \mathbb{C},\ n \in \mathbb{N}.\eqno(2.23)$$
When $n=1$ it reduces easily to the functional equation (1.6) for the Riemann zeta function via Entry 2.4.3.3 in [9], Vol. I.   In the meantime,   the odd case of $s= 2m+1,\ m \in \mathbb{Z}$ can be treated in a similar manner.  Indeed, we have from (2.9) 

$$  \beta (2(n+1), 2m+1) = {4^{n-2m} (-1)^{n+1} \over   (2n+1)! }\  \sum_{k=1}^{n+1}  t(2(n+1), 2k) \left( 2^{2k} - 2^{2m+1}\right)$$

$$\times  \zeta (2(1+m-k)),\  m \in \mathbb{Z},\ n \in \mathbb{N}_0.\eqno(2.24)$$ 
When $m$ is a negative   integer, i.e. $m= -r,\ r \in \mathbb{N}$, we take into account  the trivial zeros of the Riemann zeta function at even negative integers to find the values

$$ \beta (2(n+1), 1-2r) = 0,\quad \quad r \in \mathbb{N},\quad  n \in \mathbb{N}_0,$$ 
which can be counted as  trivial zeros of the generalized Dirichlet beta function (1.1).  Since $\zeta(0)= -1/2$, the case $m=0$ gives, in turn, the value 

$$  \beta (2(n+1),1) = {4^n  (-1)^{n} \  t(2(n+1), 2) \over   (2n+1)! }.$$
On the other hand, Entry 2.4.3.5 in [9], Vol. I implies

$$ \beta (2(n+1), 1) = {(2n)!!\over (2n+1)!!}.$$
This means that  $t (2(n+1), 2) = (-1)^{n} (n!)^2.$  Finally,  for positive integers $m$ we get the equality 

$$  \beta (2(n+1), 2m+1) = {2^{1+2(n-m)} (-1)^{n}  \over   (2n+1)! }\  \sum_{k=\  m- \min(m,n)}^{m}  t(2(n+1), 2(m+1-k)) \left( 1- 2^{1-2k}\right) \zeta (2k),\ m \in \mathbb{N},\ n \in \mathbb{N}_0.\eqno(2.25)$$

In the same manner we treat the odd case of the generalized Dirichlet beta function $ \beta(2n+1, s)$, writing, accordingly,   (see (1.24), (2.3), Entry 2.4.3.2 in [9], Vol. I)

$$ \beta (2n+1, s) = {2^{2n-s+1}\over \pi (2n)! }\ \cos\left({\pi s\over 2}\right) \int_0^\infty x^{-s}  \Gamma\left( n+{1\over 2} +ix\right)  \Gamma\left( n+{1\over 2}  - ix\right) dx$$

$$ = {2^{2n-s+1}\over  (2n)! }\ \cos\left({\pi s\over 2}\right) \int_0^\infty {x^{-s}\over \cosh (\pi x)} \   \prod_{k=1}^{n} \left[ \left(k- {1\over 2}\right)^2+ x^2 \right] dx$$

$$=  {2^{2n-s+1}\over  (2n)! }\ \cos\left({\pi s\over 2}\right) \sum_{k=0}^n (-1)^{n+k} t(2n+1, 2k+1)  \int_0^\infty {x^{2k -s}\over \cosh (\pi x)} \ dx $$

$$=  {2^{2n-s+2} (-1)^{n}\over  (2n)! }\ \cos\left({\pi s\over 2}\right) \sum_{k=0}^{n}  {(-1)^{k}\over \pi^{2k-s+1}}\  t(2n+1, 2k+1) \Gamma(2k-s+1) \beta (2k-s+1),\  0< {\rm Re}\ s < 1,\ n \in \mathbb{N}_0.$$
As in the even case we continue analytically this equality to  the whole complex plane, involving the functional equation for the Dirichlet beta function (1.4).  Hence it yields the identity

$$  (2n)!\  \beta(2n+1, s) = 2  (-1)^{n}\  \sum_{k=0}^{n}  \ 4^{n-k} \ t(2n+1, 2k+1) \  \beta (s-2k),\   s \in \mathbb{C},\ n \in \mathbb{N}_0.\eqno(2.26)$$ 
 If $s$ is an even integer $s= 2m,\ m \in \mathbb{Z}$, it becomes 
$$ \beta(2n+1, 2m) = {2^{2n+1} (-1)^{n}\over  (2n)! }\  \sum_{k=0}^{n}  \ 4^{-k}\  t(2n+1, 2k+1)  \beta (2(m-k)) ,\quad m \in \mathbb{Z},\ n \in \mathbb{N}_0.$$ 
Then  if $m$ is a nonpositive  integer, i.e. $m= -r,\ r \in \mathbb{N}_0$, we take into account  the corresponding values of the Dirichlet beta function  in terms of Euler numbers  [9], Vol. I

$$\beta (- 2n) = {1\over 2} \ E_{2n},\ \  \beta (- 2n-1) = 0, \quad n \in \mathbb{N}_0,$$

$$\beta ( 2n+1) = { (-1)^n E_{2n} \over 2 (2n)!} \ \left({\pi\over 2}\right)^{2n+1},$$
to derive

$$ \beta(2n+1, -2r) = { (-1)^{n}\over  (2n)! }\  \sum_{k=0}^{n} 4^{n-k} \  t(2n+1, 2k+1)  E_{2(r+k)},\quad r, n \in \mathbb{N}_0.$$ 
When, in turn,  $m \in \mathbb{N}$ we write

$$  \beta(2n+1, 2m) =  \  {2 (-1)^{n}\over (2n)!}  \  \sum_{k= m-n}^{m}  4^{n+k-m} \  t(2n+1, 2(m-k)+1)  \beta (2 k),\quad m \in \mathbb{N},\ n \in \mathbb{N}_0.$$ 
For $s=2m+1$ we have, accordingly, 

$$  \beta (2n+1, 2m+1) = {2^{2n-1} (-1)^{n}\over   (2n)! }\  \sum_{k=m-\min(m,n)}^{m}  \  t(2n+1, 2(m-k)+1) { (-1)^k E_{2k} \over 4^{k} (2k)!} \  \pi^{2k+1},\   m, n  \in \mathbb{N}_0.\eqno(2.27)$$ 
Using the orthogonality relation for the central factorial numbers and their properties [10], we have reciprocally from (2.26)

$$  \sum_{k=0}^{n} (-1)^{k}  2^{- 2k-1} (2k)!  T(2n+1, 2k+1) \beta (2k+1, s) =   \sum_{k=0}^{n}   T(2n+1, 2k+1) \sum_{q=0}^{k}  \  t(2k+1, 2q+1) 4^{-q} \beta (s-2q)$$

$$=    \sum_{q=0}^{n} 4^{-q} \beta (s-2q) \sum_{k=q}^{n}   T(2n+1, 2k+1)\  t(2k+1, 2q+1) =  \sum_{q=0}^{n} 4^{-q} \beta (s-2q) \delta_{2n+1,2q+1}$$

$$= 4^{-n} \beta (s-2n),$$ 
i.e.

$$ 2 \beta (s) =   \sum_{k=0}^{n} (-1)^{k}  4^{n- k} (2k)! \ T(2n+1, 2k+1) \beta(2k+1, s+2n), \quad  s \in \mathbb{C},\ n \in \mathbb{N}_0.$$
Hence  the functional equation for the Dirichlet beta function (1.4) implies, in turn,

$$   \left({\pi\over 2}\right)^{-s} \  \Gamma( s)  \sin\left({\pi s\over 2}\right) \beta (s) =   \sum_{k=0}^{n} (-1)^{k}  2^{2(n- k)-1} (2k)! \ T(2n+1, 2k+1) \beta (2k+1, 2n+1-s).$$
Therefore as in the even case we find the following functional equation for the generalized Dirichlet beta function (1.1)

$$ \sum_{k=0}^{n} (-1)^{k}  4^{-k} (2k)! \ T(2n+1, 2k+1)\bigg[   \left({\pi\over 2}\right)^{-s} \  \Gamma( s)  \sin\left({\pi s\over 2}\right)  \beta (2k+1, s+2n) \bigg.$$

$$\bigg. - \beta (2k+1, 2n+1-s) \bigg]= 0,\ \quad s \in \mathbb{C},\ n \in \mathbb{N}_0.\eqno(2.28)$$
Moreover, it implies immediately the equalities 

$$ \sum_{k=0}^{n} (-1)^{k}  4^{n- k} (2k)! \ T(2n+1, 2k+1) \beta (2k+1, 2(n+m)) =  2 \beta (2m),\ m \in \mathbb{Z},\eqno(2.29)$$

$$ \sum_{k=0}^{n} (-1)^{k}  4^{n- k} (2k)! \ T(2n+1, 2k+1) \beta(2k+1, 1+ 2(n+m)) =  { (-1)^m E_{2m} \over  (2m)!} \ \left({\pi\over 2}\right)^{2m+1},\   m \in \mathbb{N}_0,$$

$$ \sum_{k=0}^{n} (-1)^{k}  4^{- k} (2k)! \ T(2n+1, 2k+1) \beta(2k+1, 1+ 2(n-m)) = 0,  \   m \in \mathbb{N}.$$
In particular, for the Catalan constant $(m=1)$ we have the representation

$$ G= \sum_{k=0}^{n} (-1)^{k}  2^{2(n- k)-1} (2k)! \ T(2n+1, 2k+1) \beta (2k+1, 2(n+1)),\quad n \in \mathbb{N}_0. $$

\section{Properties of  $  \zeta(z,s)$}

In the same manner as in (2.1) we deduce the recurrence relation for the generalized Riemann zeta function (1.2). Precisely, it has

$$\zeta (z,s) =   z^2  \zeta (z, s+2) +  z(z+1) \  \zeta (z+2, s+2), \quad  {\rm Re}\ s >   {\rm Re}\ z > 0. \eqno(3.1)$$

{\bf Theorem 2}. { \it  Let $n \in \mathbb{N}_0$ and $c_{n,k}(\alpha)$ be defined in $(2.6)$. Then the following identity holds}

$$  \sum_{k=0}^n  {(-1)^k\over 2^k} \  \left({z\over 2}\right)_{k}\ \ c_{n,k} \left({z\over 2}\right) \  \zeta(z+2k,s+2n) =  4^{-n} \  \zeta(z,s),\ {\rm Re}\ s  >  {\rm Re}\  z > 0.\eqno(3.2)$$

\begin{proof}  In fact,  assuming  the condition $1 >  {\rm Re}\ s > {\rm Re}\ z > 0$ in (1.2), we  employ as in (2.3) the Parseval equalities for the Fourier cosine and sine transforms and the slightly corrected Entry 2.5.46.15 in [9], Vol. I. Hence we find from (1.2) the representations

$$ \zeta(z,s) = {2^{z-s}\  \cos\left(\pi s/2\right)  \over  \pi \ \Gamma(z) \  \cos\left(\pi z/2\right)}\ \int_0^\infty x^{-s}  \cosh(\pi x)\  \Gamma\left( {z\over 2}+ix\right)  \Gamma\left( {z\over 2} - ix\right) dx$$

$$ = {2^{z-s} \  \sin\left(\pi s/2\right)  \over  \pi \ \Gamma(z) \  \sin\left(\pi z/2\right)}\ \int_0^\infty x^{-s}  \sinh(\pi x)\  \Gamma\left( {z\over 2}+ix\right)  \Gamma\left( {z\over 2} - ix\right) dx.\eqno(3.3)$$
Then, subtracting two equalities,  we have 

$$  2^{s+1-z}  \pi\ \Gamma (z) \  {\sin\left(\pi (s-z) /2 \right)\over  \sin\left(\pi s \right)}\  \zeta(z,s) =  \int_0^\infty x^{-s}  e^{-\pi x} \  \Gamma\left( {z\over 2}+ix\right)  \Gamma\left( {z\over 2} - ix\right) dx,\eqno(3.4)$$
and in the same manner as in (2.4) with the use of the duplication formula for the gamma function  it gives the following integral representation

$$  2^{s- z/2}  \pi   \  \Gamma \left({z\over 2} \right)\ { \sin\left(\pi (s-z) /2 \right)\over \sin\left(\pi s \right)}\  \zeta(z,s)  =  \int_0^\infty x^{-s}  e^{-\pi x} \   \int_0^\infty t^{z/2-1} e^{-t} K_{ix}(t) dt dx,\eqno(3.5)$$
where \  $0  <  {\rm Re}\ z < {\rm Re}\ s < 1.$ Therefore,  letting again $z+2k$ instead of $z$ in (3.5),   we  arrive at the companion equality to  (2.5)

$$ 2^{-z/2} \pi\  \Gamma  \left({z\over 2}\right) {\sin\left(\pi (s-z) /2 \right)\over  \sin\left(\pi s \right)} \sum_{k=0}^n  {(-1)^k\over 2^k} \   \left({z\over 2}\right)_{k}  \  c_{n,k}\left({z\over 2}\right)\  \zeta(z+2k,s)  $$

$$ =  2^{-s} \int_0^\infty x^{-s}  e^{-\pi x} \   \int_0^\infty t^{z/2 -1} e^{-t} p_n\left(t; \ {z\over 2}\right) K_{ix}(t) dt dx,\eqno(3.6)$$
where its left-hand side is analytic in the region $ {\rm Re}\ z > 0, \  {\rm Re}\ z+ 2n  <   {\rm Re}\ s,\ s \notin \mathbb{N}$. Hence  by virtue of the  integration by parts in the integral by $t$ on the right-hand side of (3.6), the use of the properties for the sequence $\{ p_n(t;\alpha)\}_{n\ge 0}$ and a  comparison with (3.5) under assumptions $s\notin \mathbb{N},\  s-z \notin  2 \mathbb{N}$ or $s\in \mathbb{N},\  s-z \in  2 \mathbb{N}$  it yields

$$ \sum_{k=0}^n  {(-1)^k\over 2^k} \  \left({z\over 2}\right)_{k}\ \ c_{n,k} \left({z\over 2}\right) \  \zeta(z+2k,s)   =  4^{-n} \  \zeta(z,s-2n)$$
which  implies  (3.2) after a simple substitution and can be continued analytically  to the whole region ${\rm Re}\ s  >  {\rm Re}\  z > 0$. 

\end{proof} 

Letting $z=2(n+1)$ in (3.4),  we are doing similar to (2.8), employing (1.22) and Entry 2.4.10.13 in [9], Vol. I  to deduce

$$  \zeta(2(n+1),s) = {2^{2n-s+2} (-1)^{n+1} \over \pi (2n+1)! }\ \cos\left({\pi s\over 2}\right) \int_0^\infty x^{-s}  e^{-\pi x} \  \Gamma\left( n+1+ix\right)  \Gamma\left(n+1 - ix\right) dx$$

$$=  {2^{2n-s+2} (-1)^{n+1} \over (2n+1)! }\ \cos\left({\pi s\over 2}\right)  \int_0^\infty {x^{1-s} e^{-\pi x} \over \sinh(\pi x)} \   \prod_{k=1}^{n} \left[ k^2+ x^2\right] dx$$

$$=  {2^{2n-s+2}   \over (2n+1)! }\ \cos\left({\pi s\over 2}\right)  \sum_{k=0}^n (-1)^{k+1}\  t (2(n+1), 2(k+1)) \int_0^\infty {x^{1+2k-s} e^{-\pi x} \over \sinh(\pi x)} \  dx $$

$$=  {2^{2n+3}  \over (2n+1)! }\ \cos\left({\pi s\over 2}\right)  \sum_{k=1}^{n+1}  {(-1)^{k}\ 4^{- k} \over \pi^{2k-s} }  t (2(n+1), 2k)  \Gamma (2k-s) \zeta( 2k-s), $$
i.e. we establish the identity

$$  \zeta(2(n+1),s) = {2^{2n+3}  \over (2n+1)! }\ \cos\left({\pi s\over 2}\right)  \sum_{k=1}^{n+1}  {(-1)^{k} \over4^k \  \pi^{2k-s} } \  t (2(n+1), 2k)  \Gamma (2k-s) \zeta( 2k-s),\ n \in \mathbb{N}_0.\eqno(3.7) $$
Moreover, involving the functional equation (1.6), it gives

$$  \zeta(2n,s) =  {2^{2n-s}  \over (2n-1)! } \sum_{k=1}^{n}   t (2n, 2k)\  \zeta( 1+s-2k),\  s \in \mathbb{C} ,\  n \in \mathbb{N},\eqno(3.8) $$
where $\zeta(2n,s)$ has simple poles at $s= 2m,\ m= 1,2,\dots, n$ with residues 

$${4^{n-m}\  t (2n, 2m) \over (2n-1)! },$$
and, accordingly, for instance,

$$  \zeta(2n, - 1-2m) = 0,\quad  m \in \mathbb{N}_0,\  n \in \mathbb{N},$$

$$  \zeta(2n,  2m - 1) = {2^{2n-2m+1}  \over (2n-1)! } \sum_{k=1}^{m}   t (2n, 2k)\  \zeta( 2(m-k)).\eqno(3.9)$$
Comparing with (2.9), we establish another identity

$$ \beta(2n,s) + (-1)^{n}  \zeta(2n,s) =  { (-1)^{n} \over   (2n-1)! }\  \sum_{k=1}^{n} 4^{n+k-s}\  t(2n, 2k)  \zeta (1+s-2k), \ s \in \mathbb{C},\ n \in \mathbb{N}.\eqno(3.10)$$
Then, taking the corresponding formulae from the previous section one gets the corresponding results for the generalized zeta function (1.2).  Employing the functional equation (2.23),  we find 

$$ \sum_{k=1}^{n}  \ 4^{-k}  (2k-1)! \  T(2n, 2k)\bigg[   \pi^{-s}  \left( 1-2^{s}\right) \ \cos\left({\pi s\over 2}\right) \Gamma(s)   \zeta (2k,  2n+s-1)\bigg.$$

$$\bigg. - 2^{-s}   \left( 1-2^{1-s}\right) \  \zeta (2k,  2n-s)\bigg] $$

$$=  \sum_{k=1}^{n}   T(2n, 2k)\bigg[  \pi^{-s}  \left( 1-2^{s}\right) \ \cos\left({\pi s\over 2}\right) \Gamma(s)   \sum_{q=1}^{k} 4^{q-s+1-2n}\  t(2k, 2q)  \zeta (2n+s-2q)\bigg.$$

$$\bigg.  - 2^{-s}   \left( 1-2^{1-s}\right) \  \sum_{q=1}^{k} 4^{q+s-2n}\  t(2k, 2q)   \zeta (1+2n-s-2q)\bigg].$$
Interchanging the order of summation on the right-hand side of the latter equality and using orthogonality relations for the central factorial numbers [10], the parity properties (see above) and functional equation (1.6), it becomes

$$ \sum_{k=1}^{n}  \ 4^{-k}  (2k-1)! \  T(2n, 2k)\bigg[   \pi^{-s}  \left( 1-2^{s}\right) \ \cos\left({\pi s\over 2}\right) \Gamma(s)   \zeta (2k,  2n+s-1)\bigg.$$

$$\bigg. - 2^{-s}   \left( 1-2^{1-s}\right) \  \zeta (2k,  2n-s)\bigg] $$

$$=   \pi^{-s}  \left( 1-2^{s}\right) \ \cos\left({\pi s\over 2}\right) \Gamma(s)  \sum_{q=1}^{n} 4^{q-s+1-2n}\  \zeta (2n+s-2q)  \sum_{k=q}^{n}   T(2n, 2k)  t(2k, 2q) $$

$$\bigg.  - 2^{-s}   \left( 1-2^{1-s}\right) \  \sum_{q=1}^{n} 4^{q+s-2n} \zeta (1+2n-s-2q) \sum_{k=q}^{n}   T(2n, 2k)  t(2k, 2q)$$

$$=    \pi^{-s} 2^{-2s+2-2n}  \left( 1-2^{s}\right) \ \cos\left({\pi s\over 2}\right) \Gamma(s) \  \zeta (s)  -  2^{s- 2n}  \left( 1-2^{1-s}\right)  \zeta (1-s)$$

$$= 2^{s-2n}  \left( 2^{1-2s}-1\right)  \zeta (1-s).$$
Thus we derive a companion identity  to (2.23) for the generalized zeta function (1.2) 

$$ \sum_{k=1}^{n}  \ 4^{-k}  (2k-1)! \  T(2n, 2k)\bigg[   \pi^{-s}  \left( 1-2^{s}\right) \ \cos\left({\pi s\over 2}\right) \Gamma(s)   \zeta (2k,  2n+s-1)\bigg.$$

$$\bigg. - 2^{-s}   \left( 1-2^{1-s}\right) \  \zeta (2k,  2n-s)\bigg] = 2^{s-2n}  \left( 2^{1-2s}-1\right)  \zeta (1-s),\quad s \in \mathbb{C},\ n \in \mathbb{N}.\eqno(3.11)$$
The  case $s=0$ reduces to  the equality (see (3.8), (3.9))

$$ \sum_{k=1}^{n-1}  \ 4^{n-k}  (2k-1)! \  T(2n, 2k) \ \bigg[ \log 2\   \zeta (2k,  2n-1) -   \zeta (2k,  2n)\bigg] $$

$$=  \sum_{k=1}^{n-1}   t(2n, 2k)  \zeta (1+ 2(n-k))  - \log 2 \bigg[ 1+ 2 \sum_{k=1}^{n}   t (2n, 2k)\  \zeta( 2(n-k))\bigg] ,\quad  n \in \mathbb{N}.\eqno(3.12)$$
Since via (3.8), (1.6) 

$$  \pi^{-s}  \left( 1-2^{s}\right) \ \cos\left({\pi s\over 2}\right) \Gamma(s)  \sum_{k=1}^{n}  \ 4^{-k}  (2k-1)! \  T(2n, 2k)  \zeta (2k,  2n+s-1) $$

$$=  (2\pi)^{-s} 2^{-2n+1}  \left( 1-2^{s}\right) \ \cos\left({\pi s\over 2}\right) \Gamma(s)  \sum_{k=1}^{n}   T(2n, 2k)  \sum_{q=1}^{k}   t (2k, 2q)\  \zeta( s+ 2(n-q)) $$

$$=  (2\pi)^{-s}  2^{-2n+1} \left( 1-2^{s}\right) \ \cos\left({\pi s\over 2}\right) \Gamma(s) \zeta( s) = 4^{-n}  \left( 1-2^{s}\right)\zeta(1-s),$$
we let  $s= - 2m,\ m \in \mathbb{N}$ in (3.11) to obtain  the following representation for the odd zeta-values (compare with (2.29)) 

$$  \zeta (2m+1) =    \sum_{k=1}^{n}  \ 4^{n+m-k}  (2k-1)! \  T(2n, 2k)  \  \zeta (2k,  2(m+n)).\eqno(3.13)$$

A companion formula for (2.28) one obtains, letting in (3.4) $z= 2n+1$ and employing (1.24).   Hence as above

$$  {2^{s- 2n-1} (-1)^{n+1}  \pi\  (2n)! \over \sin\left(\pi s /2 \right)} \  \zeta(2n+1,s) $$

$$=  \int_0^\infty x^{-s}  e^{-\pi x} \  \Gamma\left( {1\over 2}+n +ix\right)  \Gamma\left( {1\over 2} +n - ix\right) dx$$

$$=  \int_0^\infty {x^{-s} e^{-\pi x}  \over \cosh (\pi x)} \   \prod_{k=1}^{n} \left[ \left(k- {1\over 2}\right)^2+ x^2 \right] dx$$

$$=  \sum_{k=0}^n (-1)^{n+k} t(2n+1, 2k+1)  \int_0^\infty {x^{2k -s} e^{-\pi x}  \over \cosh (\pi x)} \ dx. $$
Therefore,  appealing to (1.6) and Entry 2.4.10.17 in [9], Vol. I, we find

$$ \zeta(2n+1,s) =  {2^{2n+1} \over (2n)! }   \sum_{k=0}^n  t(2n+1, 2k+1) \bigg(  2^{-2k} - 2^{-s}\bigg) \zeta (s-2k),\  s  \in \mathbb{C},  n \in \mathbb{N}_0.\eqno(3.14)$$
Similarly, $ \zeta(2n+1,s)$ has simple poles at $s= 2m+1,\ m= 0, 1,\dots, n$ with residues

$$  {4^{n-m} \over (2n)! }\  t(2n+1, 2m+1).$$
In the meantime, reciprocally, from (3.14) 

$$ 2^{1-s}  \bigg(  2^{s}- 1 \bigg) \zeta (s)  = \sum_{k=0}^n 4^{n-k} \  T(2n+1, 2k+1) \ (2k)! \ \zeta(2k+1,s+2n),\eqno(3.15)$$
and again via (1.6) and a simple change of variable  

$$   2  \bigg(  2^{1-s}- 1 \bigg) \pi^{-s} \cos \left({\pi s\over 2}\right) \Gamma(s) \zeta (s) $$

$$ = \sum_{k=0}^n 4^{n-k} \  T(2n+1, 2k+1) \ (2k)! \ \zeta(2k+1,1-s+2n).\eqno(3.16)$$
Consequently, it is not difficult to  derive the following functional equation

$$\sum_{k=0}^n 4^{n-k} \  T(2n+1, 2k+1) \ (2k)! \bigg[  \pi^{-s} \  \Gamma(s)  \cos \left({\pi s\over 2}\right)\zeta(2k+1,s+2n) \bigg. $$

$$\bigg. +  \zeta(2k+1,1-s+2n) \bigg] = \zeta(1-s),\quad   s \in \mathbb{C},\  n \in \mathbb{N}_0.\eqno(3.17)$$
Taking $s= -2m,\ m \in \mathbb{N}$ in (3.17) and using the fact (by virtue of (1.6), (3.15)) that

$$\sum_{k=0}^n 4^{n-k} \  T(2n+1, 2k+1) \ (2k)! \lim_{s\to -2m}  \Gamma(s)    \zeta(2k+1,s+2n)  = 2^{1+2m}  \bigg(  2^{-2m}- 1 \bigg) \lim_{s\to -2m}  \Gamma(s)   \zeta(s) $$

$$=  (-1)^m \bigg(  2^{-2m}- 1 \bigg)  \pi^{-2m} \zeta (2m+1),$$ 
the odd zeta values can be written in the form (a companion of (3.13))

$$  \bigg( 2^{2m+1} -1 \bigg)  \zeta(2m+1) =  \sum_{k=0}^n 4^{n+m-k}  T(2n+1, 2k+1) \ (2k)!   \zeta(2k+1, 1+ 2(n+m)), \ m \in \mathbb{N},\  n \in \mathbb{N}_0.\eqno(3.18)$$
Finally in this section, a comparison of (2.21), (2.29), (3.13) and (3.18) suggest the following equalities, involving the odd zeta and even beta values for $m \in \mathbb{N},\  n \in \mathbb{N}_0$

$$2 \beta (2m) \pm \bigg( 2^{2m+1} -1 \bigg)  \zeta(2m+1) = \sum_{k=0}^n 4^{n-k}  T(2n+1, 2k+1) \ (2k)! $$

$$\times  \bigg[ (-1)^{k} \beta (2k+1, 2(n+m)) \pm  4^m \zeta(2k+1, 1+ 2(n+m)) \bigg],\eqno(3.19)$$

$$4^{-m} \ \zeta (2m+1) = \sum_{k=1}^{n}  \ 4^{m+n-k}  (2k-1)! \  T(2n, 2k)  \  \bigg[  \zeta (2k,  2(m+n))  + (-1)^k \beta (2k, 2(n+m)) \bigg],\eqno(3.20)$$

$$\left( 2-4^{-m}\right) \ \zeta (2m+1) =   \sum_{k=1}^{n}   4^{m+n-k}  (2k-1)!\  T(2n, 2k)\ \bigg[  \zeta (2k,  2(m+n)) -  (-1)^k \beta (2k, 2(n+m))\bigg],\eqno(3.21)$$

$$  \sum_{k=1}^{n}  \ 4^{n+m-k}  (2k-1)! \  T(2n, 2k)  \bigg[   \left( 1-4^{-m}\right)  \zeta (2k,  2(m+n)) + (-1)^k  \beta (2k, 2(n+m))\bigg] = 0.\eqno(3.22)$$

\section{Beta- and zeta-values at integers} 

In this section we will find new relations and representations for the beta- and zeta-values at integers, involving, in particular, the Kontorovich-Lebedev transform (1.11) and integrals with the algebraic and hyperbolic functions. Denoting polynomials (1.7) with $\alpha=0$ as $p_n(x; 0) \equiv p_n(x)$, we employ the integral  from [12]

$${\tau^{2n-1}\over \sinh\pi\tau} = {(-1)^n\over \pi}\int_0^\infty e^{-x}K_{i\tau}(x)p_n(x){dx\over
x}, \quad  \tau >0,\ n \in \mathbb{N}.\eqno(4.1)$$
Hence we integrate by $\tau$ over $\mathbb{R}_+$ in (4.1),  and after  simple substitutions we obtain the equalities

$$\int_0^\infty {\tau^{2n-1}\over \cosh(\pi\tau)} d\tau = {(-1)^n  \over \pi\  4^{n-1}}\int_0^\infty \sinh(\pi\tau) \int_0^\infty e^{-x}K_{2i\tau}(x)p_n(x){dx d\tau \over x}, \quad  n \in \mathbb{N},\eqno(4.2)$$

$$\int_0^\infty {\tau^{2n}\over \cosh(\pi\tau)} d\tau = {(-1)^n  \over \pi\  4^{n-1}}\int_0^\infty \tau \sinh(\pi\tau) \int_0^\infty e^{-x}K_{2i\tau}(x)p_n(x){dx d\tau \over x}, \quad  n \in \mathbb{N}.\eqno(4.3)$$
Hence,  combining  with (1.1), we  establish  iterated integrals for the beta-values at even and odd integers

$$\beta(2n) =  {(-1)^n (\pi/2)^{2n-1} \over (2n-1)!}\int_0^\infty \sinh(\pi\tau) \int_0^\infty e^{-x}K_{2i\tau}(x)p_n(x){dx d\tau \over
x}, \quad  n \in \mathbb{N},\eqno(4.4)$$

$$\beta(2n+1)  = {2 (-1)^n  (\pi/2)^{2n} \over (2n)! } \int_0^\infty \tau \sinh(\pi\tau) \int_0^\infty e^{-x}K_{2i\tau}(x)p_n(x){dx d\tau \over x}, \quad  n \in \mathbb{N}_0.\eqno(4.5)$$

{\bf Theorem 3.} {\it  Iterated integrals $(4.4), (4.5)$ can be written in terms of the series, involving Riemann's zeta-values }

$$\beta(2n) =   { 1 \over  2}  \sum_{k=n}^\infty  \binom{2k}{ 2n-1}  {2^{2k+1}-1  \over 16^{k}}  \zeta(2k+1),\quad   n \in \mathbb{N},\eqno(4.6)$$

$$\beta(2n+1) =   2  \sum_{k=n+1}^\infty  \binom{2k-1}{ 2n}  {4^{k}-1  \over 16^{k}}  \zeta(2k),\quad   n \in \mathbb{N}_0.\eqno(4.7)$$

\begin{proof}  Indeed, writing the series for the hyperbolic sine, we get (4.4) in the form

$$\beta(2n) =  {(-1)^n \pi^{2n-1} \over  2^{2n} (2n-1)!} \lim_{N\to \infty} \int_0^\infty \sum_{k=0}^N {(\pi \tau)^{2k+1}\over 2^{2k+1} (2k+1)!} \int_0^\infty e^{-x}K_{i\tau}(x)p_n(x){dx d\tau \over x},\eqno(4.8)$$
where the passage to the limit under the integral can be  motivated by the dominated convergence, involving the inequality (1.100) in [13] for the modified Bessel function.   Then  we  appeal to (1.7), (1.8), (1.9)  to write  the $n$-th iterate ${\mathcal A}^{n} \ K_{i\tau}(x)= \tau^{2n} K_{i\tau}(x), \ n \in \mathbb{N}_0.$  Consequently, integrating by parts in the integral by $x$ in (4.8), eliminating the integrated terms,  we derive 

$$\beta(2n) =  {(-1)^n  \over  2 (2n-1)!}   \sum_{k=0}^\infty  {\pi^{2(n+k)}\over 4^{n+k} (2k+1)!} \int_0^\infty \tau  \int_0^\infty e^{-x} 
{\mathcal A}^{k} K_{i\tau}(x)p_n(x){dx d\tau \over x}$$

$$ =  {(-1)^n  \over  2 (2n-1)!}   \sum_{k=0}^\infty  {\pi^{2(n+k)}\over 4^{n+k} (2k+1)!} \int_0^\infty   \int_0^\infty  \tau K_{i\tau}(x)  {\mathcal A}^{k} \left( e^{-x}  p_n(x) \right) {dx d\tau \over x}$$

$$ =  { 1 \over  2 (2n-1)!}   \sum_{k=0}^\infty  {\pi^{2(n+k)} (-1)^{n+k} \over 4^{n+k} (2k+1)!} \int_0^\infty   \int_0^\infty  \tau K_{i\tau}(x)  e^{-x}  p_{k+n} (x)  {dx d\tau \over x}.$$
The latter double integral is calculated in [12] in terms of the values of the Riemann zeta function (1.5), and we have the formula 

$$(-1)^n  (2n)!  \  \left(2^{2n+1} - 1\right)\   {\zeta(2n+1)\over (2 \pi)^{2n} } =  \int_0^\infty \int_0^\infty \tau K_{i\tau} (x) e^{-x} p_n(x){d\tau dx\over  x}, \ n \in \mathbb{N}.\eqno(4.9)$$
Thus it leads immediately to the identity (4.6). Analogously one  proves  (4.7), appealing to the identity [12] 

$$\frac{2^{2n}  - 1}{2^{2(n-1)}}\  (-1)^n  (2n-1)! { \zeta(2n)\over \pi^{2n}}  =  \int_0^\infty e^{-2x} p_n(x){dx\over x}, \ n \in \mathbb{N}\eqno(4.10)$$
and the  integral representation for polynomials $p_n$ in terms of the Kontorovich-Lebedev transform

$$p_n(x)= {2 (-1)^n\over \pi} e^x \int_0^\infty \tau^{2n} K_{i\tau}(x) d\tau,\quad x >0.\eqno(4.11)$$

\end{proof} 

On the other hand, the dominated convergence permits to write (4.4)  in the form 

$$\beta(2n) =  {(-1)^n (\pi/2)^{2n} \over 2 (2n-1)!} \lim_{\alpha \to \pi/2-} {1\over \alpha} \int_0^\infty \sinh\left(\alpha \tau \right) \int_0^\infty e^{-x}K_{i\tau}(x)p_n(x){dx d\tau \over x}$$

$$=  {(-1)^n (\pi/2)^{2n} \over 2 (2n-1)!} \lim_{\alpha \to \pi/2- } {1\over \alpha} \bigg[ \int_0^\infty \cosh\left(\alpha \tau \right) \int_0^\infty e^{-x}K_{i\tau}(x)p_n(x){dx d\tau \over
x}\bigg.$$

$$\bigg. - \int_0^\infty e^{-\alpha\tau}  \int_0^\infty e^{-x}K_{i\tau}(x)p_n(x){dx d\tau \over x}\bigg].$$
Then, one can interchange  the order of integration in iterated  integrals on the right-hand side of the latter equality, appealing  to  Entries 2.16.48.5, 2.16.48.19 in [9],  Vol. II and passing to the limit under the integral sign  to  obtain

$$  (-1)^n\ 2^{2n+1}  (2n-1)! {\beta(2n) \over \pi^{2n}} =  \lim_{\alpha \to \pi/2- } {1\over \alpha} \bigg[ {\pi\over 2}  \int_0^\infty e^{- 2 x \cos^2(\alpha/2) } p_n(x){dx \over  x}\bigg.$$

$$\bigg. - \alpha   \int_0^\infty {1\over \alpha^2+ t^2} \int_0^\infty e^{- 2 x \cosh^2(t/2)}\ p_n(x){dx d t  \over x}\bigg]$$

$$=   \int_0^\infty e^{-  x } p_n(x){dx \over  x} -  \int_0^\infty \int_0^\infty {e^{- 2 x \cosh^2(t/2)}\ p_n(x)\over x ((\pi/2)^2+ t^2) }\ dx d t,\ n \in \mathbb{N}.\eqno(4.12)$$
But (see [12], formula (2.35))

$$ \int_0^\infty e^{-  x } p_n(x){dx \over  x} = {B_{2n}\over 2n} \bigg[4^n- 4^{2n}\bigg],\ n \in \mathbb{N}.\eqno(4.13)$$
Moreover, denoting coefficients of the polynomials $p_n$ (see  (2.6)) by $c_{n,k}(0)\equiv c_{n,k}$,  recalling the Parseval equality for the Fourier  cosine transform, Entry 2.5.46.6 in [9], Vol. I and (1.22), we deduce
$$  \int_0^\infty \int_0^\infty {e^{- 2 x \cosh^2(t/2)}\ p_n(x)\over x ((\pi/2)^2+ t^2) }\ dx d t$$

$$=  \sum_{k= 1}^n c_{n,k} \int_0^\infty {1\over (\pi/2)^2+ t^2 } \int_0^\infty e^{- 2x \cosh^2(t/2)} x^{k-1}  dx dt $$

$$= 2 \sum_{k= 1}^n {c_{n,k} (k-1)! \over 2^k }  \int_{0}^\infty {dt \over \left( (\pi/2)^2+ 4t^2 \right) \cosh^{2k}  t}$$

$$ =  {1\over \pi}\sum_{k= 1}^n {c_{n,k}\  (k-1)! \  2^{k} \over (2k-1)! }  \int_{0}^\infty e^{-(\pi/2) y}  \Gamma\left(k+ i y\right) \Gamma\left(k- i y \right) dy$$

$$= \sum_{k= 1}^n {c_{n,k}\  (-1)^k \  (k-1)! \  2^{k} \over (2k-1)! } \sum_{q=0}^{k-1} (-1)^{q-1} t(2k, 2(q+1))   \int_{0}^\infty {e^{-(\pi/2) y}\  y^{2q+1}\over \sinh(\pi y)}\ dy.\eqno(4.14)$$
Hence Entry 2.4.10.13 in [9], Vol. I and interchange of the order of summation give the value of the latter  double integral  in terms of the Hurwitz   zeta function

$$   \int_0^\infty \int_0^\infty {e^{- 2 x \cosh^2(t/2)}\ p_n(x)\over x ((\pi/2)^2+ t^2) }\ dx d t$$

$$ = \sum_{q=0}^{n-1} (-1)^{q}  \  {(2q+1)!\over 2^{2q} \pi^{2q+2} } \bigg( \sum_{k= q}^{n-1} {c_{n,k+1}\  (-1)^k \  k! \  2^{k} \over (2k+1)! }\  t(2(k+1), 2(q+1)) \bigg) \ \zeta\left( 2(q+1),\  {3\over 4} \right),$$
and, in turn, the following representation for  the even beta-values (cf. (4.12), (4.13))

$$ \beta(2n) =   \bigg( 1-4^{-n}\bigg)\ \zeta(2n)  + {1 \over 2(2n-1)!} \sum_{q=1}^{n} (-1)^{n+q+1}  \  {(2q-1)!\over 4^{n+q-1} } $$

$$\times \bigg( \sum_{k= q}^{n} {c_{n,k}\  (-1)^k \  (k-1)! \  2^{k-1} \over (2k-1)! }\  t(2k, 2 q) \bigg) \ \zeta\left( 2 q,\  {3\over 4} \right) \ \pi^{2(n-q)},\quad  n \in \mathbb{N},\eqno(4.14)$$
where the Euler formula for even zeta values [1] is used.  In particular, for the Catalan constant ($n=1$) we have

$$G={1\over 8} \left(  \pi^2 -   \zeta\left( 2,\  {3\over 4} \right)\right)$$
which is a known result.  Now, appealing to the relation of the Hurwitz function with derivatives of the digamma function 

$$(2q-1)! \zeta\left( 2 q,\  {3\over 4} \right) = \psi^{(2q-1)}\left({3\over 4}\right),\quad q \in \mathbb{N}\eqno(4.15)$$
and K{\" o}lbig's formula  [6]

$$\psi^{(2q-1)}\left({3\over 4}\right) = {4^{2q-1}\over 2q} \bigg[ (-1)^{q-1} \pi^{2q} B_{2q} \left(4^q-1\right) - 2(2q)! \beta(2q)\bigg],\eqno(4.16)$$
equality (4.14) becomes

$$ \beta(2n) =   \bigg( 1-4^{-n}\bigg)\ \zeta(2n)  + {1 \over (2n-1)!} \sum_{q=1}^{n}   \  {(-1)^{n+q+1} \over 4^{n+1-q} \ q} $$

$$\times \bigg( \sum_{k= q}^{n} {c_{n,k}\  (-1)^k  \over (2k-1)!! }\  t(2k, 2 q) \bigg) \  \ \pi^{2(n-q)} \bigg[ (-1)^{q-1} \pi^{2q} B_{2q} \left(4^q-1\right) - 2(2q)! \beta(2q)\bigg]$$

$$=   \bigg( 1-4^{-n}\bigg)\ \zeta(2n)  + {(-1)^n \pi^{2n} \over (2n-1)!} \sum_{q=1}^{n}   \  {B_{2q} \left(4^q-1\right) \over 4^{n+1-q} \ q}  \sum_{k= q}^{n} {c_{n,k}\  (-1)^k  \over (2k-1)! !}\  t(2k, 2 q) $$

$$+  {1 \over (2n-1)!} \sum_{q=1}^{n}   \ (-1)^{n+q} (2q-1)!  \left({\pi\over 2}\right)^{2(n-q)}\ \beta(2q)  \sum_{k= q}^{n} {c_{n,k}\  (-1)^k  \over (2k-1)!! }\  t(2k, 2 q), $$
and it means the identity

$$ \sum_{q=1}^{n-1}  {4^{q-1} \over q} (-1)^{q+1} \pi^{-2q}  \sum_{k= q}^{n} {c_{n,k}\  (-1)^k  \over (2k-1)! !}\  t(2k, 2 q) \bigg[  (-1)^{q-1} \pi^{2q} B_{2q} \left(4^q-1\right)  - \  2  (2q)!  \ \beta(2q) \bigg] =0,$$
and via (4.16) it has

$$ \sum_{q=1}^{n-1} (-1)^q  (2\pi)^{2q}  \psi^{(2(n-q)-1)}\left({3\over 4}\right) \sum_{k= 0}^{q} { (-1)^{k}\ c_{n,n-k}  \over (2(n-k)-1)! !}\  t(2(n-k), 2 (n-q)) =0.\eqno(4.17)$$
It would be nice to have this combination of polygamma functions with nonzero coefficients. However,  it will not happen since 

$$\sum_{k= 0}^{q} {  (-1)^{k}\ c_{n,n-k}   \over (2(n-k)-1)! !}\  t(2(n-k), 2 (n-q)) =0,\quad  q= 1,2,\dots, n-1,\quad n \in \mathbb{N},\eqno(4.18)$$
and we suggest the interested reader to prove it directly. 

Further, returning to (2.1), it is not difficult to write the equality

$${\Gamma(s) \ \left[ (2n-1) \beta (2n-1, s) -  2n \ \beta (2n+1, s)  \right]\over  s-1} =   \int_0^\infty {x^{s-2} \sinh x \over \cosh^{2n} x} dx,\ s \neq 1,\ n \in \mathbb{N}.\eqno(4.19)$$
We are interested to take the limit when $s \to 1$ in (4.19), observing via Entry 2.4.3.5 in [9], Vol. I that the numerator on its left-hand side is zero when $s=1$.  Hence it gives the value of the following integral

$$\int_0^\infty { \sinh x \over  x \cosh^{2n} x} dx = - \gamma \left[ (2n-1) \beta (2n-1, 1) -  2n \ \beta (2n+1, 1)  \right] + (2n-1) \beta^\prime (2n-1, 1) -  2n \ \beta^\prime (2n+1, 1),$$
where $\gamma$ is Euler-Mascheroni constant and $\prime$ denotes the derivative with respect to $s$.  But (see Entry 2.4.3.5) )

$$ \beta (2n+1, 1) = {\pi\over 2} \ {(2n-1)!!\over (2n)!!}.\eqno(4.20)$$
Therefore 

$$\int_0^\infty { \sinh x \over  x \cosh^{2n} x} dx = (2n-1)\  \beta^\prime (2n-1, 1) -  2n  \beta^\prime (2n+1, 1),\quad n \in \mathbb{N}.\eqno(4.21)$$
On the other hand,  by virtue of Entries 2.4.18.4 in [9], Vol. I and   2.16.2.2 in [9], Vol. II  we derive 

$$\int_0^\infty { \sinh x \over  x \cosh^{2n} x} dx =  {1\over (2n-1)!} \int_0^\infty  { \sinh x \over  x } \int_0^\infty e^{-y \cosh x} y^{2n-1} dy dx $$

$$= {1\over (2n-1)!} \int_0^1 \int_0^\infty K_t(y) y^{2n-1} dy =  {2^{2n-2} \over (2n-1)!} \int_0^1 \Gamma\left( n+ {t\over 2}\right) \Gamma\left( n- {t\over 2}\right) dt,$$ 
where the interchange of the order of integration is allowed by the absolute convergence.  Then,  recalling (1.22), it has

$$\int_0^\infty { \sinh x \over  x \cosh^{2n} x} dx = {2^{2n-1} (-1)^{n+1} \over (2n-1)!} \sum_{k=1}^{n}   \pi^{1-2k}  t(2n, 2k ) \int_0^{\pi/2}   {t^{2k-1} \over \sin t} \ dt.$$
So, the right-hand side of the latter equality involves the moment integrals (cf. [12])  with the cosecant function which are calculated by formulae

$$ \int_0^{\pi/2}  { t^{2k-1} \over \sin t} \ dt =  2  (2k-1)! \sum_{m=1}^{k}  { (-1)^{m+1} \over (2(k-m))!}   \left({\pi\over 2}\right)^{2(k-m)} \beta( 2m),\eqno(4.22)$$

$$ \int_0^{\pi/2}  { t^{2k} \over \sin t} \ dt = 2 (-1)^k (2k)!  \left(1- 2^{-2k-1}\right) \zeta(2k+1) $$

$$+  2  (2k)! \sum_{m=1}^{k}  { (-1)^{m+1} \over (2(k-m)+1)!}   \left({\pi\over 2}\right)^{2(k-m)+1} \beta( 2m).\eqno(4.23)$$
Hence

$$ \int_0^\infty { \sinh x \over  x \cosh^{2n} x} dx =  {2^{2n}  (-1)^{n}  \over (2n-1)!} \sum_{k=1}^{n}   \pi^{1-2k} \   t(2n, 2k)   (2k-1)! \sum_{m=1}^{k}  { (-1)^{m} \over (2(k-m))!}   \left({\pi\over 2}\right)^{2(k-m)} \beta( 2m)$$

$$= {  (-1)^{n} 4^n  \over (2n-1)!} \sum_{m=1}^{n}  (-1)^m  \left(\sum_{k=0}^{n-m} { t(2n, 2(k+m) )\over 4^k (2k)!}\   \ (2(k+m)-1)!  \right)  {\beta( 2m)\over  \pi^{2m-1}},$$
i.e. we get the result

$$  \int_0^\infty { \sinh x \over  x \cosh^{2n} x} dx =  {  (-1)^{n} 4^{n}  \over (2n-1)!} \sum_{m=1}^{n}  (-1)^m   \left(\sum_{k=0}^{n-m} {  t(2n, 2(k+m) ) \over   4^k (2k)!}\   \ (2(k+m)-1)!  \right)  {\beta( 2m)\over  \pi^{2m-1}},\ n \in \mathbb{N}.\eqno(4.24)$$
For instance, $n=1$ implies the value

$$  \int_0^\infty { \sinh x \over  x \cosh^{2} x} dx  =  {4 G\over \pi}.$$
Meanwhile, recalling (2.26) and recurrence relation (1.17) for the central factorial numbers, we have  

$$  \beta(2n+1, s) =  {2  (-1)^{n} \over (2n)!}\ 4^n t(2n+1, 1)\  \beta (s) +  {2  (-1)^{n}\over (2n)!} \  \sum_{k=1}^{n}  \ 4^{n-k} \  t(2n-1, 2k-1)  \beta (s-2k) $$

$$ +  {  (-1)^{n+1} (2n-1) \over n  (2n-2)!} \  \sum_{k=1}^{n}  \ 4^{n-1-k} \ t(2n-1, 2k+1)\  \beta (s-2k)$$

$$= {2  (-1)^{n} \over (2n)!}\ 4^n  \left[t(2n+1, 1) + {1\over 4} \ (2n-1)^2 t(2n-1, 1) \right] \  \beta (s)  +  {2  (-1)^{n}\over (2n)!} \  \sum_{k=1}^{n}  \ 4^{n-k} \  t(2n-1, 2k-1)  \beta (s-2k) $$

$$+ {2n-1 \over 2n} \beta(2n-1, s).$$
Hence

$${2n-1 \over 2n} \beta(2n-1, s) -  \beta(2n+1, s) =  {2  (-1)^{n+1} \over (2n)!}\ 4^n  \left[t(2n+1, 1) + {1\over 4} \ (2n-1)^2 t(2n-1, 1) \right] \  \beta (s) $$

$$ +  {2  (-1)^{n+1}\over (2n)!} \  \sum_{k=1}^{n}  \ 4^{n-k} \  t(2n-1, 2k-1)  \beta (s-2k).$$
Moreover, appealing to (2.27) and (4.20),  we find the value for $t(2n+1, 1)$. Precisely, it has 

$$ t(2n+1, 1) = (-1)^n 4^{-n}\  [(2n-1)!!]^2.$$
Therefore we establish the following important equality 

$$(2n-1)  \beta(2n-1, s) - 2n\   \beta(2n+1, s) =  {2  (-1)^{n+1}\over (2n-1)!} \  \sum_{k=1}^{n}  \ 4^{n-k} \  t(2n-1, 2k-1)  \beta (s-2k),\quad n \in \mathbb{N}.\eqno(4.25)$$
Differentiating through in (4.25) with respect to $s$ at $s=1$ and comparing with (4.21), (4.24), we obtain the identity, involving  the even beta values and  derivatives of the Dirichlet beta function at odd negative integers

$$   \sum_{k=1}^{n}  \ 2^{1-2k} \  t(2n-1, 2k-1)  \beta^\prime (1-2k) $$

$$=  \sum_{m=1}^{n} (-1)^{m+1}  \left(\sum_{k=0}^{n-m} {t(2n, 2(k+m) ) \over 4^k (2k)!} \ (2(k+m)-1)!  \right)  {\beta( 2m)\over  \pi^{2m-1}},\ n \in \mathbb{N}.\eqno(4.26)$$
Let $n=1$ in (4.26). Then it gives the known value $\beta^\prime (-1)= 2G /\pi.$  But the known general formula for $\beta^\prime (1-2n),\ n \in \mathbb{N}$

$$\beta^\prime (1-2n) = (-1)^{n+1} (2n-1)! \left({2\over \pi}\right)^{2n-1} \beta(2n)$$
 which follows from functional equation (1.4) for the Dirichlet beta function, taking logarithm  and using the differentiation.    Thus we get
 
 $$   \sum_{k=1}^{n} (-1)^{k+1}   t(2n-1, 2k-1)  (2k-1)!  {\beta(2k) \over \pi^{2k-1}} =  \sum_{m=1}^{n} (-1)^{m+1} \sum_{k=0}^{n-m} {t(2n, 2(k+m) ) \over 4^k (2k)!} \ (2(k+m)-1)!  {\beta( 2m)\over  \pi^{2m-1}}$$
and conjecture a possibly new identity for even and odd central factorial numbers of the first kind 

$$ 4^{-k} t(2n-1, 2k-1)  =  \sum_{q=k}^{n} 4^{-q}  \binom{ 2q -1} {2k-1} \   t(2n, 2q),\quad k, n \in \mathbb{N}.\eqno(4.27)$$
 Furthermore,  (4.21) becomes 
 
 $$  h_{2n} \equiv \int_0^\infty { \sinh x \over  x \cosh^{2n} x} dx =  {  (-1)^{n} 4^n  \over (2n-1)!}  \sum_{k=1}^{n}  (-1)^k \ (2k-1)! \   t(2n-1, 2k-1)\  {\beta(2k)\over \pi^{2k-1}},\quad n \in \mathbb{N}.\eqno(4.28)$$
For $n=2$ we have the value of the integral

$$  h_4= \int_0^\infty { \sinh x \over  x \cosh^{4} x} dx =  {16\over \pi}  \left[ {\beta(4)\over \pi^2}  +  {G\over 24 } \right].$$
A reciprocal identity  to (4.28) is established by the following theorem. 

{\bf Theorem 4}.  {\it Let $n \in \mathbb{N}$ and $h_{2n}$ be defined by $(4.28)$. Then even beta-values have the  representation 

$$ {\beta(2n)\over \pi^{2n-1} } =    { 4^{-n}  \over (2n-1)!} \sum_{k=1}^{n}  d_{n,k} \   h_{2k},\eqno(4.29)$$
where }

$$d_{n,k} =  {(-1)^{n}\over 4^{k-1}} \sum_{r=1}^{k} (-1)^r \binom{2k-1}{k-r} (2r-1)^{2n-1} = \lim_{t\to 0 }\  {d^{2n-1}\over dt^{2n-1} }  \sin^{2k-1} t.\eqno(4.30)$$

\begin{proof} In fact, recalling equalities (4.6), (4.9), we have 
$$ \beta(2n)=   { 1 \over  2 (2n-1)!}   \sum_{k=n}^\infty  {\pi^{2k} (-1)^{k} \over 4^{k} (2(k-n)+1)!} \int_0^\infty   \int_0^\infty  \tau K_{i\tau}(x)  e^{-x}  p_{k} (x)  {dx d\tau \over x}$$

$$=  { 1 \over  2 (2n-1)!} \lim_{t\to \pi/2 - }   \sum_{k=n}^\infty  {t^{2k} (-1)^{k} \over (2(k-n)+1)!} \int_0^\infty   \int_0^\infty  \tau K_{i\tau}(x)  e^{-x}  p_{k} (x)  {dx d\tau \over x}$$

$$ = { 1 \over  2 (2n-1)!}   \lim_{t\to \pi/2 - }  t^{2n-1} {d^{2n-1}\over dt^{2n-1} }\sum_{k=1}^\infty  {t^{2k} (-1)^{k} \over (2k)!} \int_0^\infty   \int_0^\infty  \tau K_{i\tau}(x)  e^{-x}  p_{k} (x)  {dx d\tau \over x},\eqno(4.31)$$
where the second equality is due to the dominated convergence and,  as  will be shown below, the repeated differentiation in the open disk $|t| < \pi/2$ is  allowed due to the absolute and uniform convergence.  Precisely, taking the integral representation for polynomials $p_n$ (cf. [12])

$$p_n(x) = {2(-1)^n\over \pi}\ e^x \int_0^\infty \tau^{2n} K_{i\tau}(x) d\tau,\eqno(4.32)$$
we write with the use of Entry 2.16.48.19 in [9], Vol. II

$$\sum_{k=0}^\infty  {t^{2k}  p_k(x) \over (2k)!}=  {2\over \pi} \ e^x \sum_{k=0}^\infty  {t^{2k} (-1)^{k} \over (2k)!} \int_0^\infty \tau^{2k} K_{i\tau}(x) d\tau$$

$$= {2\over \pi} \ e^x \ \int_0^\infty \cos(t\tau)  K_{i\tau}(x) d\tau = e^{-2x\sinh^2(t/2) },\quad x >0,\quad  \left|t \right| < {\pi\over 2}\eqno(4.33)$$
which is the generating function for these polynomials,  and the interchange of the order of integration and summation is permitted by virtue of the estimate 

 $$\sum_{k=0}^\infty  {|t|^{2k}  \over (2k)!} \int_0^\infty \tau^{2k} \left| K_{i\tau}(x) \right| d\tau =  \int_0^\infty \cosh(|t|\tau) \left| K_{i\tau}(x) \right| d\tau$$

$$\le K_0(x\cos\delta)  \int_0^\infty e^{\tau( |t|- \delta)} d\tau,\quad |t| < \delta < {\pi\over 2},$$
where the inequality (1.100) from [13] is employed.  Hence, taking into account the recurrence relation [12] for $p_n$

$${p_n(x)\over x} = - \sum_{k=0}^{n-1} \binom{2n-1}{2k} p_k(x),\quad n \in \mathbb{N},$$
 one justifies the third equality in (4.31) because 

$$\sum_{k=1}^\infty  {|t|^{2k-m} \over (2k-m)!}  \int_0^\infty   \int_0^\infty  \tau \left| K_{i\tau}(x)\right|   e^{-x} | p_{k} (x) | {dx d\tau \over x}$$

$$\le {2\over \pi} \sum_{k=1}^\infty  {|t|^{2k-m} \over (2k-m)!}  \sum_{q=0}^{k-1} \binom{2k-1}{2q} \int_0^\infty  K_0(x\cos\delta) \ dx  \int_0^\infty  \tau^{2q+1}  e^{-\delta\tau} \  d\tau $$

$$<  {2\over \pi} \sum_{k=1}^\infty  {|t|^{2k-m} \over (2k-m)!}  \sum_{q=0}^{2k-1} \binom{2k-1}{q} \int_0^\infty  K_0(x\cos\delta) \ dx  \int_0^\infty  \tau^{q+1}  e^{-\delta\tau} \  d\tau $$

$$=  {2\over \pi}  \int_0^\infty  K_0(x\cos\delta) \ dx\  \sum_{k=1}^\infty  {|t|^{2k-m} \over (2k-m)!}   \int_0^\infty  \tau\  (1+\tau)^{2k-1}  e^{-\delta\tau} \  d\tau $$

$$<   {2\over \pi}  \int_0^\infty  K_0(x\cos\delta) \ dx\  \sum_{k=1}^\infty  {|t|^{2k-m} \over (2k-m)!}   \int_0^\infty    (1+\tau)^{2k}  e^{-\delta\tau} \  d\tau$$

$$ <   {2 e^{\delta_0} \over \pi}\   \int_0^\infty  K_0\left(x\cos\delta_0\right) \ dx   \int_0^\infty    (1+\tau)^{m}  e^{\tau(\delta_0 -\delta)} \  d\tau < \infty,\ |t| \le \delta_0 <  \delta < {\pi\over 2},\quad  m=1,2,\dots, 2n-1.$$
Consequently,  taking into account (4.33), we find from (4.31) the equality

$$  \beta(2n) = { 1 \over  2 (2n-1)!}   \lim_{t\to \pi/2 - }  t^{2n-1} {d^{2n-1}\over dt^{2n-1} } \int_0^\infty   \int_0^\infty  \tau K_{i\tau}(x)   \left( e^{-x \cos t} -e^{-x} \right)  {dx d\tau \over x}.\eqno(4.34)$$
In the meantime, we will fulfill  the repeated differentiation under the integral sign,  invoking Hoppe's formula which gives, in turn,

$${d^{2n-1}\over dt^{2n-1} } \left[ e^{-x \cos t} \right] =  e^{-x \cos t} \sum_{k=0}^{2n-1} {  x^k \over k!} \sum_{j=0}^k \binom{k}{j} (-1)^j \left( \cos t\right)^{k-j}{d^{2n-1}\over dt^{2n-1} }  \cos^j t.$$
Then, accordingly,

$$ \lim_{t\to \pi/2 - }  t^{2n-1} {d^{2n-1}\over dt^{2n-1} } \left[ e^{-x \cos t} \right] =  \lim_{t\to \pi/2 - }  t^{2n-1} e^{-x \cos t} \sum_{k=0}^{2n-1} {  x^k \over k!} \sum_{j=0}^k \binom{k}{j} (-1)^j \left(  \cos t\right)^{k-j} {d^{2n-1}\over dt^{2n-1} }  \cos^j t$$

$$=  \left({\pi\over 2} \right)^{2n-1} \lim_{t\to \pi/2 - } {d^{2n-1}\over dt^{2n-1} }   \sum_{k=0}^{2n-1} { (-1)^k  x^k \over k!} \ \cos^k t = 
 - \left({\pi\over 2} \right)^{2n-1} \lim_{t\to \pi/2 - } {d^{2n-1}\over dt^{2n-1} }   \sum_{k=0}^{n-1} {  x^{2k+1} \cos^{2k+1} t\over (2k+1)!},$$
and, employing   Entry 2.16.2.2 in [9], Vol. II,  we find from (4.34)

$$ {\beta(2n)\over \pi^{2n-1}} =   { 1 \over 2 (2n-1)!}   \sum_{k=1}^{n} {4^{k-1-n}  d_{n,k}\over (2k-1)!} \int_0^\infty  \tau\  \Gamma\left( k -1+ {1+i\tau\over 2} \right) \Gamma\left( k-1 + {1-i\tau\over 2}\right)  d\tau,\eqno(4.35)$$
where coefficients $d_{n,k}$ are defined in (4.30).  Meanwhile, Entry 2.5.46.5 in [9], Vol. I and integration by parts imply  the equalities

$$ {1\over (2k-1)!} \int_0^\infty  \tau \ \Gamma\left( k -1+ {1+i\tau\over 2} \right) \Gamma\left( k -1+ {1-i\tau\over 2}\right)  d\tau =  2^{3-2k}  \int_0^\infty d\tau \int_0^\infty {\sin ( \tau u) \sinh u  \over \cosh^{2k} u} du$$

$$= 2^{3-2k}  \lim_{N\to \infty}  \int_0^N  d\tau \int_0^\infty {\sin ( \tau u) \sinh u  \over \cosh^{2k} u} du = 2^{3-2k}  \lim_{N\to \infty}  \int_0^\infty {(1- \cos ( N u) ) \sinh u  \over  u \cosh^{2k} u} du$$

$$= 2^{3-2k}  \int_0^\infty {\sinh u  \over u \cosh^{2k} u} du,$$ 
where the interchange of the order of integration and passage to the limit under the integral sign are due to the dominated convergence and the Riemann-Lebesgue lemma. 
Consequently, substituting the latter expression in (4.35), we arrive at (4.29), completing the proof of Theorem 4. 

\end{proof}

A similar result is valid for the odd Riemann zeta-values.   Precisely, it has 

{\bf Theorem 5}.  {\it Let $n \in \mathbb{N}$ and $h_{2n+1}$ be defined accordingly (cf.  $(4.28)$). Then for the odd Riemann zeta-values  the  following reciprocal equalities hold  

$$  (2n)!  \  \left(2^{2n+1} - 1\right)\   {\zeta(2n+1)\over  (2\pi)^{2n} } =  (-1)^n \sum_{k=1}^{n} 2^{-k}  k!\ c_{n,k}\  h_{2k+1},\eqno(4.36)$$

$$ h_{2n+1} =  {  (-1)^n 4^{n}  \over (2n)!}  \sum_{k=1}^{n}  (-1)^k   \  t(2n, 2k) \ (2k)!  \left( 2^{2k+1} -1\right) {\zeta(2k+1)\over  (2\pi)^{2k}},\eqno(4.37)$$
where $c_{n,k}$ are coefficients of the polynomials $p_n$.}

\begin{proof} Indeed, in the same manner and justifications we have owing to  Entry 2.4.3.1 in [9], Vol. I 

 $$(-1)^n  (2n)!  \  \left(2^{2n+1} - 1\right)\   {\zeta(2n+1)\over  \pi^{2n} } = \int_0^\infty \tau \cosh \left({\pi\tau\over 2}\right) \int_0^\infty e^{-x}K_{i\tau}(x)p_n(x){dx d\tau \over x}$$

$$= (-1)^n \sum_{k=n}^\infty  {\pi^{2(k-n)} (-1)^{k} \over 4^{k-n} (2(k-n))!} \int_0^\infty \int_0^\infty \tau  e^{-x}K_{i\tau}(x)p_k(x){dx d\tau \over x}$$

$$=   (-1)^n \lim_{t\to \pi/2 - }   {d^{2n}\over dt^{2n} } \int_0^\infty   \int_0^\infty  \tau K_{i\tau}(x)   \left( e^{-x \cos t} -e^{-x} \right)  {dx d\tau \over x}.$$
Therefore  we deduce

$$ (2n)!  \  \left(2^{2n+1} - 1\right)\   {\zeta(2n+1)\over  \pi^{2n} } = \lim_{t\to \pi/2 - }   {d^{2n}\over dt^{2n} } \int_0^\infty   \int_0^\infty  \tau K_{i\tau}(x)   \left( e^{-x \cos t} -e^{-x} \right)  {dx d\tau \over x}$$

$$=   \lim_{t\to \pi/2 - }   \sum_{k=0}^{n}  {2^{2k-2} \over (2k)!} \ {d^{2n}\over dt^{2n} }  \cos^{2k} t   \int_0^\infty  \tau \Gamma\left( k + {i\tau\over 2} \right) \Gamma\left( k - {i\tau\over 2}\right) d\tau $$

$$=   {(-1)^n\over 2} \sum_{k=1}^{n}  \sum_{r=0}^{k-1} {(2(k-r))^{2n} (-1)^{k-r} \over r!  (2k-r)! }  \int_0^\infty  \tau \Gamma\left( k + {i\tau\over 2} \right) \Gamma\left( k - {i\tau\over 2}\right) d\tau $$

$$=   (-1)^{n}  2^{1+2n}  \sum_{k=1}^{n}  \sum_{r=1}^{k} (-1)^r  4^{-k} \ \binom{2k} {k-r}  r^{2n}   \int_0^\infty { \sinh u\over u  \cosh^{2k+1} u} du.$$
Hence it yields

$$ 2  (2n)!  \  \left(2^{2n+1} - 1\right)\   {\zeta(2n+1)\over  (2\pi)^{2n} } =  \sum_{k=1}^{n} 4^{1-k} \left( \sum_{r=1}^{k}  (-1)^{r+n} r^{2n} \binom{2k}{k+r} \right) h_{2k+1}.$$
But Entry 4.4.1.20 in [9], Vol. I and formula for the polynomial coefficients $c_{n,k}$ in  [12] imply

$$ \sum_{r=1}^{k}  (-1)^{r+n} r^{2n} \binom{2k}{k+r} = 2^{2k-2n-1} \lim_{x\to 0} {d^{2n}\over d x^{2n} } \ \sin^{ 2k} x =  2^{k-1} (-1)^n k!\ c_{n,k}.$$
Hence we get (4.36).  Conversely, recalling  Entries 2.4.18.4 in [9], Vol. I,    2.16.2.2 in [9], Vol. II and  (1.24), we derive 

$$ h_{2n+1} = \int_0^\infty { \sinh x \over  x \cosh^{2n+1} x} dx =  {1\over (2n)!} \int_0^\infty  { \sinh x \over  x } \int_0^\infty e^{-y \cosh x} y^{2n} dy dx $$

$$= {1\over (2n)!} \int_0^1 \int_0^\infty K_t(y) y^{2n} dy =  {2^{2n-1} \over (2n)!} \int_0^1 \Gamma\left( n+ {1+t\over 2}\right) \Gamma\left( n+ {1-t\over 2}\right) dt$$ 

$$=  {2^{2n-1} \over (2n)!} \int_0^1 \left(n- {t\over 2}\right) \Gamma\left( n -  {t\over 2}\right) \Gamma\left( n+ {t\over 2}\right) dt$$ 

$$=  {  (-1)^{n+1} \pi  \over 2 (2n)!}  \sum_{k=0}^{n-1}  4^{n-k-1}   t(2n, 2(k+1) ) \int_0^{1}   {(2n-t)\  t^{2k+1} \over \sin (\pi t/2)} \ dt $$ 

$$=  {  (-1)^{n+1} 4^n  \over 2 (2n)!}  \sum_{k=1}^{n}    {t(2n, 2k) \over  \pi^{2k-1}}\ \int_0^{\pi/2}  \left(2n- {2 t \over \pi} \right)\  {t^{2k-1} \over \sin t} \ dt. $$ 
Thus, invoking (4.22), (4.23) and simplifying the expression after substitution, we establish the identity 

$$ h_{2n+1} =   { 2^{2n+1} \over  (2n)!}   \sum_{m=1}^{n-1}  {(-1)^m\over 4^m} \   {\beta( 2(n-m)) \over  \pi^{2(n-m)-1}} \sum_{k=0}^{m} 4^k\  t(2n, 2(n-k))   { (2(n-k)-1)! \over  (2(m-k) +1)!} \ \left(2n (m-k)  +k \right) $$ 

$$ + {  (-1)^n 4^{n}  \over  (2n)!}  \sum_{k=1}^{n}  (-1)^k   \  t(2n, 2k) \ (2k)!  \left( 2^{2k+1} -1\right) {\zeta(2k+1)\over  (2\pi)^{2k}}.\eqno(4.38)$$
On the other hand, in the same manner we obtain a companion of the equality (4.21) which yields 

$$h_{2n+1} =  2n\  \beta^\prime (2n, 1) -  (2n+1)  \beta^\prime (2(n+1), 1),\quad n \in \mathbb{N}.\eqno(4.39)$$
Then,  employing  (2.9) and (1.17), we write 

$$ 2n \beta (2n, s) - (2n+1) \beta (2(n+1), s) = {4^{n-s} (-1)^{n} 2n \over   (2n-1)! }\  \sum_{k=1}^{n}  t(2n, 2k) \left( 2^{2k} - 2^{s}\right) \zeta (1+s-2k)$$

$$+ {4^{n+1-s} (-1)^{n} \over   (2n)! }\  \sum_{k=1}^{n+1}  t(2(n+1), 2k) \left( 2^{2k} - 2^{s}\right) \zeta (1+s-2k)$$

$$ = {4^{n-s} (-1)^{n} 2n \over   (2n-1)! }\  \sum_{k=1}^{n}  t(2n, 2k) \left( 2^{2k} - 2^{s}\right) \zeta (1+s-2k)$$

$$+ {4^{n+1-s} (-1)^{n} \over   (2n)! }\  \sum_{k=1}^{n+1}  t(2n, 2(k-1)) \left( 2^{2k} - 2^{s}\right) \zeta (1+s-2k)$$

$$-  {4^{n-s} (-1)^{n} 2n \over   (2n-1)! }\  \sum_{k=1}^{n+1}  t(2n, 2k) \left( 2^{2k} - 2^{s}\right) \zeta (1+s-2k)$$

$$=  {4^{n+1-s} (-1)^{n} \over   (2n)! }\  \sum_{k=1}^{n+1}  t(2n, 2(k-1)) \left( 2^{2k} - 2^{s}\right) \zeta (1+s-2k),$$
i.e.

$$  2n \beta (2n, s) - (2n+1) \beta (2(n+1), s)  = {4^{n+1-s} (-1)^{n} \over   (2n)! }\  \sum_{k=0}^{n}  t(2n, 2k) \left( 2^{2(k+1)} - 2^{s}\right) \zeta (s-2k-1).$$
Therefore, making a differentiation with respect to $s$ at $s=1$ through the latter equality and taking the known values 

$$\zeta^\prime (-2n) = {(-1)^n\over 2 (2\pi)^{2n}} \ (2n)! \zeta (2n+1),\quad n \in \mathbb{N},$$ 
we find from (4.39)

$$h_{2n+1} = {4^{n} (-1)^{n}  \over   (2n)! }\  \sum_{k=1}^{n} (-1)^k \  t(2n, 2k)\ (2k)!\  \left( 2^{2k+1} - 1 \right) {\zeta (2k+1)\over  (2\pi)^{2k}}.$$
Comparing with (4.38), it implies (4.37) because  

$$\sum_{m=1}^{n-1}  {(-1)^m\over 4^m} \   {\beta( 2(n-m)) \over  \pi^{2(n-m)-1}} \sum_{k=0}^{m} 4^k\  t(2n, 2(n-k))   { (2(n-k)-1)! \over  (2(m-k) +1)!} \ \left(2n (m-k)  +k \right) = 0.\eqno(4.40)$$ 
 Theorem 5 is proved. 
 \end{proof}
 
 As an immediate consequence of (4.40) we conjecture the following identity for central factorial numbers  
 
$$ \sum_{k=0}^{m} 4^k\  t(2n, 2(n-k)) \ \binom{2(n-k)-1}{2(n-m-1)} \ \left(2n (m-k)  +k \right) = 0,\quad  m \in \mathbb{N}_0,\ m \le n-1,\ m \in \mathbb{N}.\eqno(4.41)$$ 
Letting $n=1,2$ in (4.37),  we have the values

$$\int_0^\infty { \sinh u\over u  \cosh^{3} u} du = {7\over \pi^2} \zeta(3),$$

$$\int_0^\infty { \sinh u\over u  \cosh^{5} u} du =  {1\over 3 \pi^2} \bigg[  7\  \zeta(3) +  93 \ {\zeta(5)\over \pi^{2}}\bigg].$$

\section{ Catalan's constant}

As we could see in Section 1,  the Catalan constant $G=\beta(2)$ is defined by the series (1.3), and it is the value of the generalized Dirichlet beta function $\beta(1,2)/2$.  Our goal in this section is to find integral representations for the function $\beta(z,2)$ which will contain the corresponding expressions for the Catalan constant. In fact, taking series (1.15), we write

$$ \beta(z,2) =  2^z \sum_{n=0}^\infty {(-1)^n\  (z)_n \over n! \ (2n+z)^2} =  2^z \sum_{n=0}^\infty {(-1)^n \over n! }\ (z)_n   \int_0^1 \int_0^1 (xy)^{2n+z-1} dx dy$$

$$ = 2 \int_0^1\int_0^1 { (2xy)^{z-1}  dxdy\over (1+ x^2y^2)^{z}},\quad 0 < {\rm Re\  z} <  1,$$
where the interchange of the order of summation and integration is  justifies with the aid of the uniform convergence.  But the latter equality  holds for all $ {\rm Re\  z}  > 0$ via analytic continuation, and we obtain the identity

$$ \beta(z,2) =   2 \int_0^1\int_0^1 { (2xy)^{z-1} dxdy\over (1+ x^2y^2)^{z}},\quad  {\rm Re\  z}  > 0.\eqno(5.1)$$

Another double integral representation for $\beta(z,2)$ is given by 

{\bf Theorem 6}. {\it Let ${\rm Re\  z}  > 0$.  It has the identity

$$\beta(z,2)  =  \int_0^1\int_0^1 { (2xy)^{z-1} dxdy\over (x+y)^{z}(1-xy)}.\eqno(5.2)$$
In particular, the Catalan constant can be  represented by  the following Beukers type integral}

$$G=  {1\over 2} \int_0^1\int_0^1 { dxdy \over (x+y)(1-xy)}.\eqno(5.3)$$

\begin{proof}  In fact, the proof is based on the Fourier cosine integral for the product of gamma functions which is a reciprocal relation to  Entry 2.5.46.6 in [9], Vol. I

$${1\over \cosh^{z} x} = {2^{z-1 }\over \pi \Gamma(z)} \int_0^\infty \Gamma\left({z+i\tau\over 2}\right) \Gamma\left({z-i\tau\over 2}\right) \cos(x\tau ) d\tau.\eqno(5.4)$$ 
Then, substituting the right-hand side of (5.4) in the integral (1.1) for $\beta(z,2)$, we write the chain of equalities

$$ \beta(z,2) = {2^{z-1 }\over \pi \Gamma(z)} \lim_{N\to \infty } \int_0^N x  \int_0^\infty \Gamma\left({z+i\tau\over 2}\right) \Gamma\left({z-i\tau\over 2}\right) \cos(x\tau ) d\tau dx $$
$$= {2^{z-1 }\over \pi \Gamma(z)} \lim_{N\to \infty } \int_0^N  \int_0^\infty \Gamma\left({z+i\tau\over 2}\right) \Gamma\left({z-i\tau\over 2}\right) {d\over d\tau}  \left[ \sin(x\tau )\right]  d\tau dx $$

$$= {2^{z-1 }\over \pi \Gamma(z)} \lim_{N\to \infty }   \int_0^\infty \Gamma\left({z+i\tau\over 2}\right) \Gamma\left({z-i\tau\over 2}\right) {d\over d\tau}  \left[ {1-\cos( N \tau)\over \tau }\right]  d\tau,$$
where the interchange of the order of integration and differentiation are  justified by the absolute and uniform convergence.  Then, integrating by parts in the latter integral by $\tau$ and appealing to the Riemann-Lebesgue lemma, we obtain 

$$\beta(z,2)  =  {2^{z}\over \pi i \ \Gamma(z)}   \int_0^\infty  \Gamma\left({z+i\tau\over 2}\right) \Gamma\left({z-i\tau\over 2}\right) $$

$$\times \bigg[ \psi \left({z+i\tau\over 2}\right) - \psi \left({z- i\tau\over 2}\right)\bigg] {d\tau\over \tau},\eqno(5.5)$$
where  the difference of the digamma functions in (5.5) has the familiar series representation

$$\psi \left({z+i\tau\over 2}\right) - \psi \left({z- i\tau\over 2}\right) = \sum_{k=0}^\infty { i\tau \over ( k+ (z- i\tau)/2 ) ( k+ (z+ i\tau)/2 )}.\eqno(5.6)$$
Hence equality (5.5) becomes

$$\beta(z,2) =  {2^{z-2}\over \pi \ \Gamma(z)}   \int_0^\infty  \Gamma\left({z+i\tau\over 2}\right) \Gamma\left({z-i\tau\over 2}\right)$$

$$\times \sum_{k=0}^\infty { 1 \over ( k+ (z- i\tau)/2 ) ( k+ (z+ i\tau)/2 )} d\tau$$

$$= {2^{z-3}\over \pi \ \Gamma(z)}   \int_{-\infty}^\infty  \Gamma\left({z+i\tau\over 2}\right) \Gamma\left({z-i\tau\over 2}\right)   \int_0^1\int_0^1 { (xy)^{(z- 2)/2}  \over 1-xy}  \left({x\over y}\right)^{-i\tau /2} dxdyd\tau.\eqno(5.7)$$
 Finally, we appeal to Fubini's theorem to interchange the order of integration and we calculate the Melllin-Barnes integral by $\tau$ via Entry 8.4.2.5 in [9], Vol. III to get from (5.7) the representation
 
 $$ \beta(z,2) = 2^{z-1}   \int_0^1\int_0^1  \left( 1+ {x\over y}\right)^{ - z} {x^{z-1} y^{- 1}  \over 1-xy}   dxdy$$ 
which yields (5.2). Letting $z=1$, we find (5.3).

\end{proof}

The following theorem gives double integral relationships, involving an arbitrary polynomial $P_n$ of degree at most $n \in \mathbb{N}_0$ with rational coefficients and three constants $G, \ \zeta(2), \ \log 2$.

{\bf Theorem 7.}  {\it Let $P_n,\ n \in \mathbb{N}_0$ be an arbitrary polynomial of degree at most $n$ with rational coefficients and nonzero free term, i.e.

$$P_n(x)= \sum_{k=0}^n  a_{n,k} x^k,\quad  a_{n,k} \in \mathbb{Q}, \quad  k= 0,1,\dots,n,\quad a_{n,0} \neq 0.\eqno(5.8)$$
Then }

$$  \int_0^1\int_0^1 P_n\left( { 2xy \over 1+ x^2y^2} \right) { dxdy\over  1+ x^2y^2}  $$

$$= {1\over 2}   \int_0^1\int_0^1 P_n\left( { 2xy \over x+y} \right) { dxdy\over (x+y)( 1- xy)} =  G\  \mathbb{Q} + \log 2 \ \mathbb{Q} + \mathbb{Q},\eqno(5.9) $$

$$ \int_0^1\int_0^1 {P_n\left( x+y \right) \over (x+y)( 1- xy)} dxdy =  G\  \mathbb{Q} + \zeta (2) \ \mathbb{Q} + \mathbb{Q},\eqno(5.10) $$

$$ \int_0^1\int_0^1 {P_n\left( xy \right) \over (x+y)( 1- xy)} dxdy =  G\  \mathbb{Q} + \log (2) \ \mathbb{Q} + \mathbb{Q},\eqno(5.11) $$

$$ \int_0^1\int_0^1 {P_n\left( (1-xy) (x+y)\right) \over (x+y)( 1- xy)} dxdy =  G\  \mathbb{Q}  + \mathbb{Q},\eqno(5.12) $$

$$ \int_0^1\int_0^1 {P_n\left( xy (x+y)\right) \over (x+y)( 1- xy)} dxdy =  G\  \mathbb{Q}  + \zeta (2) \mathbb{Q} + \mathbb{Q},\eqno(5.13) $$

$$ \int_0^1\int_0^1 {P_n\left( x^2y^2 \right) \over 1+ x^2y^2} dxdy =  G\  \mathbb{Q}  + \mathbb{Q},\eqno(5.14) $$

$$  \int_0^1\int_0^1 P_n\left( { (xy)^r  \over x+y} \right) { dxdy\over (x+y)( 1- xy)} =  G\  \mathbb{Q} + \log 2 \ \mathbb{Q} + \mathbb{Q},\quad  r \in \mathbb{N},\eqno(5.15) $$

$$ \int_0^1\int_0^1 {P_n\left( (xy)^r \right) \over (1+ x^2y^2)( 1- xy)} dxdy =    G \mathbb{Q}  + \zeta (2) \mathbb{Q} + \mathbb{Q},\quad  r \in \mathbb{N},\eqno(5.16) $$

$$ \int_0^1\int_0^1 {P_n\left( (1+xy)(x+y) \right) \over (x+y)( 1- x^2y^2)} dxdy =   G  \mathbb{Q}  + \zeta (2) \mathbb{Q} + \mathbb{Q}.\eqno(5.17) $$

\begin{proof} The proof follows immediately from (2.26), (5.1), (5.2), (5.3)  and elementary integrals

$$ \int_0^1 { x^n  \over x+1}\  dx = (-1)^n \left[ \log 2    + \sum_{k=1}^n {(-1)^k\over k}\right],$$

$$ \int_0^1\int_0^1 {dxdy \over x+y} = 2\log 2,$$

$$ \int_0^1\int_0^1 {dxdy \over 1+x^2y^2} =  G,$$

$$ \int_0^1\int_0^1 {dxdy \over 1- xy} =  \zeta(2),$$

$$ \int_0^1\int_0^1 {dxdy \over 1+ xy} =  {1\over 2} \zeta(2),$$

$$ \int_0^1\int_0^1 {dxdy  \over (1+ x^2y^2)( 1- xy)}  =  {G\over 2}  + {9\over 16} \zeta(2),$$

$$ \int_0^1\int_0^1 {dxdy  \over (1- x^2y^2)( x+y)}  =  G  + {3 \over 16} \zeta(2).$$

\end{proof}

Simple sufficient conditions to guarantee the irrationality of Catalan's constant is given by

{\bf Theorem 7}.  {\it $G \notin \mathbb{Q}$ if there exists a polynomial  $(5.8)$ with integer coefficients such that

$$\sum_{k=0}^n  \sum_{m=0}^k  \sum_{r=0}^k   { (-1)^{r+k} a_{n,k+1} \over (m+r+1)(m+k-r+1)} \in \ \mathbb{Z},\eqno(5.18)$$
and the double integral $(5.12)$ converges to zero when $n \to \infty$.}

\begin{proof} In fact, recalling  equality (5.12), we have via (5.3), (5.18)

$$ \int_0^1\int_0^1 {P_n\left( (1-xy) (x+y)\right) \over (x+y)( 1- xy)} dxdy =  2 a_{n,0}\  G +  \sum_{k=0}^n  \sum_{m=0}^k  \sum_{r=0}^k   { (-1)^{r+k} a_{n,k+1} \over (m+r+1)(m+k-r+1)} $$

$$= G \mathbb{Z}+ \mathbb{Z}.\eqno(5.19)$$
Hence since the double integral in  (5.19) tends to zero when $n \to \infty$, the constant $G$ cannot be rational.

\end{proof}

\section{Ramanujan-type formulae }

Our departure point in this section will be the Ramanujan formula (1.29) written in the form (cf. [1])

$$a^{-k} \sum_{n=1}^\infty {\coth(a n) \over n^{2k+1} }=  (- b)^{-k}  \sum_{n=1}^\infty {\coth(b n) \over n^{2k+1} } - 2^{2k+1} \sum_{j=0}^{k+1} {(-1)^j B_{2j} B_{2k+2-2j}\over (2j)!(2k+2-2j)!} a^{k+1-j} b^j.\eqno(6.1)$$
However  Entry 2.5.46.3 in [9], Vol. I says 

$${a n\over 2\pi } \coth \left(a n\right)= {1\over 2\pi } +  \int_0^\infty \left({\sin ( a nx)\over \sinh (\pi x)}\right)^2 \ dx.$$
Therefore via the dominated convergence  and the equality $ab= \pi^2$ it gives 

$$ a^{-k} \sum_{n=1}^\infty {\coth(a n) \over n^{2k+1} } =  a^{-k-1} \zeta(2(k+1)) + 2 \pi a^{-k-1}  \int_0^\infty {dx \over \sinh^2(\pi x)} \sum_{n=1}^\infty {\sin^2 ( a nx) \over n^{2(k+1)} } $$

$$=  a^{-k-1} \zeta(2(k+1)) -  \pi^2 a^{-k-2}  \int_0^\infty {dy \over \sinh^2(b y)} \sum_{n=1}^\infty {\cos(2\pi n y) -1 \over n^{2(k+1)} }.$$
Then, employing Euler formula for even zeta values and  using Entry 5.4.2.7 in [9], Vol. I, we have 

$$  \sum_{n=1}^\infty {\coth(a n) \over n^{2k+1} } = {\zeta(2(k+1))\over a } -   {(-1)^{k}\  2^{2k+1}\pi^{2(k+2)} \over a^2  (2(k+1))!}  \int_0^\infty {dy \over \sinh^2(b y)} \left[\ B_{2(k+1)} \left( y\right) -   B_{2(k+1)} \right]$$

$$+  {(-1)^{k}\  2^{2k+1}\pi^{2(k+2)} \over a^2  (2(k+1))!}  \int_1^\infty {dy \over \sinh^2(b y)} \left[  \ B_{2(k+1)} \left( y\right)   -   B_{2(k+1)} \left( \{y\}\right)\right],\eqno(6.2)$$
where $\{ \}$ is the fractional part of the number and $B_n(y),\ n \in \mathbb{N}$ are Bernoulli polynomials.   Hence we easily find the formula 

$$  \sum_{n=1}^\infty {\coth(a n) \over n^{2k+1} } =  {(-1)^{k}\  (2\pi)^{2(k+1)}  \over 2 a  (2(k+1))!} \bigg[ B_{2(k+1)}   +  b   \int_0^\infty {dy \over \sinh^2(b y)} \left[  \ B_{2(k+1)}    -   B_{2(k+1)} \left( \{y\}\right)\right] \bigg].\eqno(6.3)$$
Moreover, from the Lerch formula (1.35) it has for $a=b=\pi,\ k= 2n+1$

$$  {1\over 2}  \sum_{m=0}^{2n+1}  (-1)^{m}  \  {B_{2m} \ B_{4n+4-2m} \over (2m)! (4n+4-2m)!} = {  B_{4(n+1)}  \over  2 (4(n+1))!}  +  {\pi \over   (4(n+1))!}  \int_0^\infty {dy \over \sinh^2(\pi y)} \left[  \ B_{4(n+1)}    -   B_{4(n+1)} \left( \{y\}\right)\right],$$
i.e. the value of the integral 

$$ \pi \int_0^\infty {dy \over \sinh^2(\pi y)} \left[  \ B_{4(n+1)}    -   B_{4(n+1)} \left( \{y\}\right)\right] = {1\over 2}  \sum_{m=1}^{2n+1}  (-1)^{m} \binom{4(n+1)} {2m}  \  B_{2m} \ B_{4n+4-2m}.\eqno(6.4)$$
Meanwhile, employing  the explicit formula for the Bernoulli polynomials [12] and the corrected Entry 2.4.3.3 in [9], Vol. I,  we find, involving the generalized Riemann zeta function (1.2)  

$${(-1)^{k}\  2^{2k+1}\pi^{2(k+2)} \over a^2  (2(k+1))!}  \int_0^\infty {dy \over \sinh^2(b y)} \left[\ B_{2(k+1)} \left( y\right) -   B_{2(k+1)} \right]$$

  $$= {(-1)^{k}\  2^{2k+1} a^k \over b^{k+1} }  \sum_{m=0}^{2k} {B_m\  b^{m} \over m!}  \zeta(2, 2k+3-m) =  (-1)^{k+1} \left({a\over b}\right)^k  \bigg[  \zeta(2k+1)\bigg.$$
 $$\bigg. + (-1)^k (2\pi)^{2k+1}  \sum_{m=0}^{k}  (-1)^{m+1}  \  {B_{2m} \ B_{2(k-m+1)} \over (2m)! (2(k-m+1))!} \  \left({b\over \pi}\right)^{2m-1} \bigg].\eqno(6.5)$$
  Thus  we derive from (6.2) 
 
 $$  \sum_{n=1}^\infty {\coth(a n) \over n^{2k+1} } =  (-1)^{k}\    \zeta(2k+1)\  \left({\pi\over b}\right)^{2k} $$
 
 $$+ (2\pi) ^{2k+1} \ \sum_{m=0}^{k+1}  (-1)^{m+1}  \  {B_{2m} \ B_{2(k-m+1)} \over (2m)! (2(k-m+1))!} \ 
  \left({\beta\over \pi}\right)^{2(m-k)-1} $$

$$+ 2\  {(-1)^{k}\  (2\pi)^{2k} b^2  \over  (2(k+1))!}  \int_0^\infty {dy \over \sinh^2(b y)} \left[  \ B_{2(k+1)} \left( y\right)   -   B_{2(k+1)} \left( \{y\}\right)\right]$$

$$=  (-1)^{k}\    \zeta(2k+1)\  \left({\pi\over b}\right)^{2k} $$
 
 $$+ (2\pi) ^{2k+1} \ \sum_{m=0}^{k+1}  (-1)^{m+1}  \  {B_{2m} \ B_{2(k-m+1)} \over (2m)! (2(k-m+1))!} \ 
  \left({b\over \pi}\right)^{2(m-k)-1} $$

$$+ 2\  {(-1)^{k}\  (2\pi)^{2k} b^2  \over  (2k+1)!}  \int_0^\infty {dy \over \sinh^2(b y)} \int_0^y  \ B_{2k+1} \left( t\right) dt $$

$$  + 2\  {(-1)^{k}\  (2\pi)^{2k} b^2  \over  (2(k+1))!}  \int_0^\infty {dy \over \sinh^2(b y)} \left[  \ B_{2(k+1)}    -   B_{2(k+1)} \left( \{y\}\right)\right].$$
Moreover,  from (1.29) we find 

$$ \sum_{n=1}^\infty {\coth(b n) \over n^{2k+1} }  =    \zeta(2k+1)\   +   { 2^{2k+1} b^{2(k+1)}  \over  (2(k+1))!}  \int_1^\infty {dy \over \sinh^2(b y)} \left[  \ B_{2(k+1)} \left( y\right)   -   B_{2(k+1)} \left( \{y\}\right)\right].\eqno(6.6)$$
Therefore, invoking  properties of the Bernoulli polynomials, it yields  

$$ \sum_{n=1}^\infty {n^{-2k-1} \over  e^{2b n} -1 }  =     { 2^{2k} b^{2(k+1)}  \over  (2(k+1))!}  \int_1^\infty {dy \over \sinh^2(b y)} \left[  \ B_{2(k+1)} \left( y\right)   -   B_{2(k+1)} \left( \{y\}\right)\right]$$

$$=  { 2^{2k} b^{2(k+1)}  \over  (2k+1)!}  \int_1^\infty {dy \over \sinh^2(b y)} \int_{\{y\}}^y  B_{2k+1} \left( u\right) du =  { 2^{2k} b^{2(k+1)}  \over  (2k+1)!} \sum_{m=1}^\infty  \int_{m}^{m+1}  {dy \over \sinh^2(b y)} \sum_{q=1}^m \left( \{y\} + q-1\right)^{2k+1} $$

$$=  { 2^{2k} b^{2(k+1)}  \over  (2k+1)!}  \sum_{m=1}^\infty  \int_{0}^{1}  {dy \over \sinh^2(b (y+m))} $$

$$\times \sum_{q=1}^m    (y+q-1)^{2k+1} =    { 2^{2k} b^{2(k+1)}  \over  (2k+1)!}  \int_{0}^{1}  dy \sum_{q=1}^\infty  \  \sum_{m=q}^\infty    {\left( y + q-1\right)^{2k+1} \over \sinh^2(\beta (y+m))} =   \int_{0}^{1}  dy \sum_{q=0}^\infty  \  \sum_{m=0}^\infty    {\left( y + q\right)^{2k+1} \over \sinh^2(b (y+m+q+1))}$$

$$=  { 2^{2k} b^{2(k+1)}  \over  (2k+1)!}  \sum_{m=0}^\infty   \sum_{q=0}^\infty  \  \int_{q}^{q+1}     { y^{2k+1} \over \sinh^2(b (y+m+1))} \ dy=  { 2^{2k} b^{2(k+1)}  \over  (2k+1)!} \int_{0}^{\infty} \left( \sum_{m=1}^\infty    {1 \over \sinh^2(b (y+m))}\right)\   y^{2k+1} \ dy,$$
where the interchange of the summation and integration is by virtue of the dominated convergence.  Thus it gives  

$$ \sum_{n=1}^\infty {n^{-2k-1} \over  e^{2b n} -1 }  =   { 2^{2k} b^{2(k+1)}  \over  (2k+1)!}  \int_{0}^{\infty} \left( \sum_{m=1}^\infty    {1 \over \sinh^2(b (y+m))}\right)\   y^{2k+1} \ dy = { 2^{2k} b^{2(k+1)}  \over  (2k+1)!}   \sum_{m=1}^\infty m^{2(k+1)}  \int_{0}^{\infty}   {y^{2k+1}  \over \sinh^2(b\ m (y+1))} \ dy.$$
Appealing to  the Mellin-Barnes type representation for $\hbox{cosech}^2 x$, which can be obtain from the corrected Entry 2.4.3.3 in [9] via the inverse Mellin transform, we have

$$ \sum_{n=1}^\infty {n^{-2k-1} \over  e^{2bn} -1 }  =  { 2^{2k+1} b^{2(k+1)}  \over  \pi  i(2k+1)!}   \sum_{m=1}^\infty m^{2(k+1)}  \int_{0}^{\infty}   y^{2k+1}  \int_{\mu-i\infty}^{\mu+i\infty} \Gamma(s) \zeta (s-1) (2b\ m (y+1))^{-s} ds dy,\quad \mu  > 2k+3.$$
Interchanging the order of integration and summation due to the absolute and uniform convergence and calculating the elementary beta-integral, we find the known identity

$$ \sum_{n=1}^\infty {n^{-2k-1} \over  e^{2b n} -1 }  =  { 1 \over 2\pi i}  \int_{\mu - 2(k+1) -i\infty}^{\mu- 2(k+1) +i\infty}\  \Gamma(s) \zeta (s) \zeta( s + 2k+1)\ (2 b )^{-s} ds.\eqno(6.7)$$ 
 Now,  differentiating in (6.7) $2k+1$ times with respect to $b$ under the integral sign owing to the absolute and uniform convergence, we derive

$$ {d^{2k+1}\over db^{2k+1}}  \sum_{n=1}^\infty {n^{-2k-1} \over  e^{2b n} -1 }  = -  { 1 \over 2\pi i}  \int_{\mu - 2(k+1) -i\infty}^{\mu- 2(k+1) +i\infty}\  \Gamma(s+2k+1) \zeta (s) \zeta( s + 2k+1)\  2^{-s} b^{-s-2k-1} ds$$

$$= -  { 2^{2k+1}  \over 2\pi i}  \int_{\mu - 1 -i\infty}^{\mu- 1 +i\infty}\  \Gamma(s) \zeta (s) \zeta( s - 2k-1)\   (2b)^{-s} ds.$$ 
Hence the functional equation (1.6) for the Riemann zeta function implies

$$ {d^{2k+1}\over db^{2k+1}}  \sum_{n=1}^\infty {n^{-2k-1} \over  e^{2b n} -1 }  = {  (-1)^{k}  \pi^{-2(k+1)}\over 4\pi i}  \int_{\mu - 1 -i\infty}^{\mu- 1 +i\infty}\   \zeta (1-s) \zeta(2k+2- s )\Gamma(2k+2-s)   (2\alpha)^{s} ds$$

$$= {  (-1)^{k} 2^{2k+1} a^{k+1}\over 2\pi i\ b^{k+1} }  \int_{2k+3-\mu -i\infty}^{2k+3-\mu +i\infty}\   \zeta (s-2k-1) \zeta( s )\Gamma(s)   (2 a)^{-s} ds.$$
Taking $2k+3 < \mu < 2k+4$, we shift the contour to the right and minding  the residue at simple poles $s=0,\ 2(k+1)$ of the integrand, we have

$$  {d^{2k+1}\over db^{2k+1}}  \sum_{n=1}^\infty {n^{-2k-1} \over  e^{2b n} -1 }  = {  (-1)^{k} 2^{2k+1} a^{k+1}\over 2\pi i\ b^{k+1} }  \int_{\mu-1 -i\infty}^{\mu-1 +i\infty}\   \zeta (s-2k-1) \zeta( s )\Gamma(s)   (2 a)^{-s} ds$$

$$+   (-1)^{k} 2^{2k} \zeta(-1-2k)  \left({a\over b}\right)^{k+1} -  {  (-1)^{k} \over 2\  \pi^{2(k+1)} } \ (2k+1)!\  \zeta (2(k+1)),$$
i.e. after simplification it gives the formula

$$  {d^{2k+1}\over db^{2k+1}}  \sum_{n=1}^\infty {n^{-2k-1} \over  e^{2b n} -1 }  =  (-1)^{k+1}  \left({a\over b}\right)^{k+1}\ {d^{2k+1}\over d a^{2k+1}}  \sum_{n=1}^\infty {n^{-2k-1} \over  e^{2a n} -1 }  - \bigg[ 1 +  (-1)^k \left({a\over b}\right)^{k+1} \bigg] {2^{2k-1}  B_{2(k+1)}  \over k+1}.\eqno(6.8)$$
In particular, for $k= 2n,\ a = b=\pi$ we find

$$ {d^{4n+1}\over d b^{4n+1}}  \sum_{m=1}^\infty {m^{-4n-1} \over  e^{2b m} -1 }  \Big |_{b=\pi} \ = -\  {2^{4n-1}\  B_{4n+2}  \over 2n+1}.\eqno(6.9)$$ 
More generally,  taking some $m \in \mathbb{N}$, we have

$$ {d^m\over db^m}  \sum_{n=1}^\infty {n^{-2k-1} \over  e^{2b n} -1 }  =   { (-1)^m \over 2\pi i}  \int_{\mu - 2(k+1) -i\infty}^{\mu- 2(k+1) +i\infty}\  \Gamma(s+m) \zeta (s) \zeta( s + 2k+1)\  2^{-s} b^{-s-m} ds$$

$$=  { (-1)^m (2b)^{2k+1} b^{-m}  \over 2\pi i}  \int_{\mu - 1 -i\infty}^{\mu - 1 +i\infty}\  \Gamma(s+m-2k-1) \zeta (s) \zeta( s - 2k-1)\   (2 b)^{-s} ds.$$ 
Hence the functional equation (1.6) yields

$$ {d^m\over db^m}  \sum_{n=1}^\infty {n^{-2k-1} \over  e^{2b n} -1 }  = { (-1)^{k} b^{2k-m+1} \pi^{-2k-2}  \over 4\pi i} \int_{\mu - 1 -i\infty}^{\mu- 1 +i\infty}\   {\Gamma(1-s)\over \Gamma(2+2k-m-s)} \zeta (1-s) \zeta(2k+2- s )\Gamma(2k+2-s)   (2 a)^{s} ds$$

$$= { (-1)^{k} 2^{2k+1} b^{k-m}  a^{k+1} \over 2\pi i}  \int_{2-\mu -i\infty}^{2-\mu +i\infty}\     {\Gamma(s) \Gamma(s+2k+1)\over \Gamma(s+2k+1-m)} \  \zeta (s) \zeta( s+2k+1 )   (2 a)^{-s-2k-1} ds$$

$$= { (-1)^{k+1} b^{k-m}  a^{k+1} \over 2\pi i} {d^{2k+1}\over d a^{2k+1} }  \int_{2-\mu -i\infty}^{2-\mu +i\infty}\     {\Gamma^2 (s) \over \Gamma(s+2k+1-m)} \  \zeta (s) \zeta( s+2k+1 )   (2a)^{-s} ds.$$
Let $m > 2k+1$. Then we write

$$ {d^m\over db^m}  \sum_{n=1}^\infty {n^{-2k-1} \over  e^{2b n} -1 }  = { (-1)^{k+1} b^{k-m}  a^{k+1} \over 2\pi i} {d^{2k+1}\over da^{2k+1} }  \int_{2-\mu -i\infty}^{2-\mu +i\infty}\    (s+2k+1-m)_{m-2k-1} \Gamma(s)  \zeta (s) \zeta( s+2k+1 )   (2 a)^{-s} ds$$

$$= { (-1)^{k+m} b^{k-m}  a^{k+1} \over 2\pi i} {d^{m}\over da^{m} }  a^{m-2k-1}  \int_{2-\mu -i\infty}^{2-\mu +i\infty}\  \Gamma(s) \zeta (s) \zeta( s+2k+1 )   (2 a)^{-s} ds.\eqno(6.10)$$
Hence we get via (1.6), (6.7) and the relation $ab = \pi^2$

$$ {d^m\over db^m}  \sum_{n=1}^\infty {n^{-2k-1} \over  e^{2b n} -1 }  =  (-1)^{m} b^{-m}  a {d^{m}\over da^{m} } a^{m-1} \sum_{n=1}^\infty {n^{-2k-1} \over  e^{2b n} -1 },\quad  m > 2k+1.\eqno(6.11)$$
Moreover,  employing the Ramanujan formula (1.29), equality (6.11) becomes 

$$ {d^m\over db^m}  \sum_{n=1}^\infty {n^{-2k-1} \over  e^{2b n} -1 }  =  (-1)^{m} b^{-m}  \alpha {d^{m}\over d a^{m} } a^{m-1} \bigg[ \bigg[ (-1)^k \pi^{2k} a^{-2k} - 1\bigg]\ {\zeta(2k+1)\over 2}  \bigg.$$

$$\bigg. +   2^{2k} \sum_{q=0}^{k+1} (-1)^{k+q}  {B_{2q} B_{2(k+1-q)} \ \pi^{2(k+q)} \over  (2q)! (2(k+1-q))!} \  a^{1-2q}+  (-1)^k \pi^{2k} a^{-2k}  \sum_{n=1}^\infty {n^{-2k-1} \over  e^{2a n} -1 }\bigg]$$

$$=   (-1)^{m} b^{-m}  \alpha {d^{m}\over da^{m} } a^{m-2k-1} \bigg[ \bigg[ (-1)^k \pi^{2k} - a^{2k} \bigg]\ {\zeta(2k+1)\over 2}  \bigg.$$

$$\bigg. +   2^{2k} \sum_{q=0}^{k+1} (-1)^{q+1}  {B_{2q} B_{2(k+1-q)} \ \pi^{2(2k+1-q)} \over  (2q)! (2(k+1-q))!} \  a^{2q-1}+  (-1)^k \pi^{2k}  \sum_{n=1}^\infty {n^{-2k-1} \over  e^{2a n} -1 }\bigg]$$

$$= (-1)^{m} \pi^{2k} b^{-m} a {d^{m}\over da^{m} } \bigg[ 2^{2k} \sum_{q=0}^{k+1} (-1)^{q+1}  {B_{2q} B_{2(k+1-q)} \ \pi^{k+1-q} \over  (2q)! (2(k+1-q))!} \  \alpha^{m-2(k-q+1)}+  (-1)^k  a^{m-2k-1} \sum_{n=1}^\infty {n^{-2k-1} \over  e^{2a n} -1 }\bigg]$$

$$= (-1)^{k+m} \pi^{2k}  b^{-m} \alpha \bigg[ 2^{2k}  {B_{2(k+1)} \ m! \over  (2(k+1))! } +  {d^{m}\over da^{m} }  a^{m-2k-1} \sum_{n=1}^\infty {n^{-2k-1} \over  e^{2a n} -1 }\bigg]$$

$$= (-1)^{k+m} \pi^{2k}  m!\ b^{-m}  \alpha \bigg[ 2^{2k}   {B_{2(k+1)} \over  (2(k+1))! } +  \sum_{q= 1}^{m-2k}  \binom{m-2k-1} {q-1} {\alpha^{q-1} \over (q+2k)!} {d^{q+2k}\over da^{q+2k}}  \sum_{n=1}^\infty {n^{-2k-1} \over  e^{2a n} -1 }\bigg].$$
Thus  we establish  the following consequence of the Ramanujan formula (1.29) for the derivatives of the involved series (1.31) 

$$ {(-1)^{k+m} \over m!} \pi^{2(m-k)} a^{-m-1}  F_{2k+1}^{(m)} (b)  -  \sum_{q= 1}^{m-2k}  \binom{m-2k-1} {q-1} {a^{q-1} \over (q+2k)!}   F_{2k+1}^{(q+2k)} (a) $$

$$=  { 4^k B_{2(k+1)} \over  (2(k+1))! }, \quad m > 2k+1.\eqno(6.12)$$
In particular, taking $a= b =\pi$, we find  from (6.11) the equality

$$\bigg[ {(-1)^{k+m} -1\over m!} \pi^{m-2k-1}   F_{2k+1}^{(m)} (a)  - \sum_{q= 1}^{m-2k-1} \binom{m-2k-1} {q-1} {\pi^{q-1}\over (q+2k)!}  F_{2k+1}^{(q+2k)} (a) \bigg] \Big |_{\alpha=\pi}=     {  4^{k}\ B_{2(k+1)} \over  (2(k+1))! }.$$
 When $m=2(k+1)$ it reads

$$\bigg[ \pi \left[ (-1)^{k} -1\right]   {d^{2(k+1)}\over d\alpha^{2(k+1)} }  \sum_{n=1}^\infty {n^{-2k-1} \over  e^{2\alpha n} -1 }  - 2(k+1)  {d^{2k+1}\over d\alpha^{2k+1}}  \sum_{n=1}^\infty {n^{-2k-1} \over  e^{2\alpha n} -1 }\bigg] \Big |_{\alpha=\pi}=  4^{k}\ B_{2(k+1)}$$
which confirms for even $k$ by (6.9).

Further, for the case $m < 2k+1,\ m \in \mathbb{N}$ we return to (6.10) to observe  in a similar manner via properties of the Mellin transform (cf. [14])

$$ {d^m\over db^m}  \sum_{n=1}^\infty {n^{-2k-1} \over  e^{2b n} -1 }  =  { (-1)^{k+1} b^{k-m}  a^{k+1} \over 2\pi i} {d^{2k+1}\over da^{2k+1} }  \int_{2-\mu -i\infty}^{2-\mu +i\infty}\     {\Gamma (s) \over (s)_{2k+1-m}} \  \zeta (s) \zeta( s+2k+1 )   (2a)^{-s} ds$$

$$= { (-1)^{k+m} b^{k-m}  a^{k+1} \over 2\pi i} {d^{m}\over da^{m} }  a^{m-2k-1}  \int_{2-\mu -i\infty}^{2-\mu +i\infty}\  \Gamma(s) \zeta (s) \zeta( s+2k+1 )   (2a)^{-s} ds$$

$$=(-1)^{m} b^{-m}  a {d^{m}\over da^{m} } a^{m-1} \sum_{n=1}^\infty {n^{-2k-1} \over  e^{2b n} -1 }$$

$$=(-1)^{m} b^{-m}  a {d^{m}\over da^{m} } a^{m-1} \bigg[ \bigg[ (-1)^k \pi^{2k} a^{-2k} - 1\bigg]\ {\zeta(2k+1)\over 2}  \bigg.$$

$$\bigg. +   2^{2k} \sum_{q=0}^{k+1} (-1)^{k+q}  {B_{2q} B_{2(k+1-q)} \ \pi^{2(k+q)} \over  (2q)! (2(k+1-q))!} \  a^{1-2q}+  (-1)^k \pi^{2k} a^{-2k}  \sum_{n=1}^\infty {n^{-2k-1} \over  e^{2a n} -1 }\bigg]$$

$$=   (-1)^{m} b^{-m}  a {d^{m}\over da^{m} } a^{m-2k-1} \bigg[  {(-1)^k \pi^{2k}\over 2}  \ \zeta(2k+1) \bigg.$$

$$\bigg. +   2^{2k} \sum_{q=0}^{k+1} (-1)^{q+1}  {B_{2q} B_{2(k+1-q)} \ \pi^{2(2k+1-q)} \over  (2q)! (2(k+1-q))!} \  a^{2q-1}+  (-1)^k \pi^{2k}  \sum_{n=1}^\infty {n^{-2k-1} \over  e^{2a n} -1 }\bigg].\eqno(6.13)$$
Therefore, fulfilling the differentiation, we deduce the Ramanujan-type formula 

$$ (-1)^k {\pi^{2(m-k)}\over m!}    {d^m\over db^m}  \sum_{n=1}^\infty {n^{-2k-1} \over  e^{2b n} -1 }  =  { 1\over 2} \binom{2k}{m}   a^{m-2k}   \ \zeta(2k+1) $$

$$+   (-1)^{m}   2^{2k}  a^{1+m}  \  {B_{2(k+1)} \over  (2(k+1))!} +   {2^{2k-1} \over m!}     a^{m+1}   \sum_{q=0}^{k} {(-1)^{q+1}\over q+1}  \ {B_{2(q+1)} B_{2(k-q)}\over (2q+1-m)!\  (2(k-q))!}  \ \left({\pi\over a}\right)^{2(q+1)} $$

$$+    \sum_{q=0}^m  {(-1)^q\over q!}  \binom{2k-q}{2k-m} \  a^{m+q-2k}\   {d^{q}\over da^{q} }  \sum_{n=1}^\infty {n^{-2k-1} \over  e^{2a n} -1 },\quad m \in \mathbb{N},\ m < 2k+1.\eqno(6.14)$$
The case $m=0$ concludes by (6.13) and (1.29), and we have 

$$  \sum_{n=1}^\infty {n^{-2k-1} \over  e^{2b n} -1 }  =  \left[ (-1)^{k} (\pi/a)^{2k}-1\right]\  {\zeta(2k+1)\over 2 }  +  (-1)^{k}  2^{2k} \pi^{2k} a \  {B_{2(k+1)}\  \over  (2(k+1))!} $$

$$+   {(-1)^{k}  2^{2k-1} \pi^{2(k+1)}\over (2k+1)!}   \sum_{q=0}^{k} (-1)^{q+1} \binom{2k+1}{2q+1}  \ B_{2(q+1)} B_{2(k-q)}  \ \pi^{2q}  \  {a^{-1- 2q}\over q+1}  $$

$$+  (-1)^{k}  \pi^{2k}  \  a^{-2k}\     \sum_{n=1}^\infty {n^{-2k-1} \over  e^{2a n} -1 }.$$
Taking $a = b=\pi$, the previous formula drives to the Lerch identity (1.35). Meanwhile, (6.14) reads in this case (see (1.31))

$$ \left( 1- (-1)^{k+m}\right) F_{2k+1}^{(m)} (\pi)  = (-1)^{k} \ m!  \sum_{q=0}^{m-1}   {(-1)^q\over q!}  \binom{2k-q}{m-q} \  \pi^{q-m}\  F_{2k+1}^{(q)} (\pi)  $$

$$+ { (-1)^{k} (2k)!\over 2 (2k-m)!} \   \pi^{-m}  \ \zeta(2k+1) +   2^{2k} \pi^{2k+1-m}  \bigg[ (-1)^{k+m} - \binom{2k+1}{m} \bigg]\ \  {B_{2(k+1)}\  m! \over  (2(k+1))!} $$

$$+   {(-1)^{k}  2^{2k-1} \pi^{2k+1-m} \over (2k+1-m)!}   \sum_{q=1}^{k} {(-1)^{q}\over q} \  \binom{2k+1-m}{2q-1-m}  \ B_{2q} B_{2(k+1-q)},\quad m \in \mathbb{N}\eqno(6.15)$$
and can be called the generalized Lerch formula. Let $m=1,  k= 2r,\ r \in \mathbb{N}.$  Then we have the equality

$$  F_{4r+1}^\prime (\pi)  = {2r\over \pi}    \sum_{n=1}^\infty {n^{-4r-1} \over  e^{2\pi n} -1 }  $$

$$+    {r\over \pi}     \zeta(4r+1) -   2^{4r-1} \pi^{4r}   \bigg[   {B_{4r+2}\   \over  (4r+2)!} +   \sum_{q=0}^{2r} (-1)^{q} \  {(2q+1)  B_{2(q+1)} B_{4r- 2q}\over (2(q+1))! (4r-2q)!}\bigg].$$
But from the Ramanujan formula (1.29) we observe that 

$$\sum_{q=0}^{2r} (-1)^{q}\  { B_{2q+2} B_{4r-2q}\over (2q+2)! (4r-2q)! }  = { B_{4r+2} \over (4r+2)!}.$$
Hence it has

$$    \zeta(4r+1) =  {2^{4r-1}\over r} \  \pi^{4r+1}   \sum_{q=0}^{2r} (-1)^{q} \  { B_{2(q+1)} B_{4r- 2q}\over (2q+1)! (4r-2q)!}-  2 \sum_{n=1}^\infty {n^{-4r-1} \over  e^{2\pi n} -1 } + {\pi\over r} \ F_{4r+1}^\prime (\pi).$$
But if we let  $m=2,  k= 2r$, we find 

$$2  \sum_{n=1}^\infty {n^{-4r-1} \over  e^{2\pi n} -1 } - {\pi\over r} \ F_{4r+1}^\prime (\pi)  +   \zeta(4r+1) +   {2^{4r} \pi^{4r+1}\over r (4r-1)}  \  {B_{4r+2} \over  (4r+2)!} $$

$$-  { 2^{4r-1} \pi^{4r+1}\over r (4r-1)}     \sum_{q=0}^{2r} (-1)^{q} \   {2q (2q+1) B_{2(q+1)} B_{4r-2q}\over (2(q+1))! (4r-2q)!} =0,$$
i.e.

$$    \zeta(4r+1) =  {2^{4r-1}\over r (4r-1)} \  \pi^{4r+1}   \sum_{q=0}^{2r} (-1)^{q} \  { (2q-1)\  B_{2(q+1)} B_{4r- 2q}\over (2q+1)! (4r-2q)!}-  2 \sum_{n=1}^\infty {n^{-4r-1} \over  e^{2\pi n} -1 }  $$

$$+ {\pi\over r} F_{4r+1}^\prime (\pi)$$
which undoubtedly confirms the previous equality  and can be  a companion of the Lerch formula for $\zeta(4r+1)$, because

   $$\sum_{q=0}^{2r-1} (-1)^{q} \  {  B_{2(q+1)} B_{4r- 2q}\over (2q+1)! (4r-2q-1)!} =0.$$
Now, letting $m=k$, (6.15) yields

$$   \binom{2k}{k}  \ \bigg[ \zeta(2k+1) +   2 \sum_{n=1}^\infty {n^{-2k-1} \over  e^{2\pi n} -1 } \bigg] $$

$$+    (-1)^k   (2\pi)^{2k+1}  \  {B_{2(k+1)} \over  (2(k+1))!} \bigg[ 1- \binom{2(k+1)}{k+1} +     {(2(k+1))!  \over 2 k! \ B_{2(k+1)}} \  \sum_{q=0}^{k-1} {(-1)^{k+q+1}\over q+1}  \ {B_{2(q+1)} B_{2(k-q)}\over (2q+1-k)!\  (2(k-q))!}  \bigg] $$

$$+    2 \sum_{q=1}^{k-1}  {(-1)^q\over q!}  \binom{2k-q}{k} \  \pi^{q}\   {d^{q}\over d\alpha^{q} }  \sum_{n=1}^\infty {n^{-2k-1} \over  e^{2\alpha n} -1 }\Big |_{\alpha=\pi} = 0.$$
For $k=1$ it confirms immediately via  the equality (1.36).

Now, returning to (6.7),  let us take the $m$-th  derivative  ($m \in \mathbb{N}$) with respect to $\alpha$. So, it  implies
$$   {d^m\over d\alpha^m}  \sum_{n=1}^\infty {n^{-2k-1} \over  e^{2\beta n} -1 }  =  { 1 \over 2\pi i}  {d^m\over d\alpha^m}   \int_{\gamma - 2(k+1) -i\infty}^{\gamma- 2(k+1) +i\infty}\  \Gamma(s) \zeta (s) \zeta( s + 2k+1)\ 2^{-s} \pi^{-2s}  \alpha^{s} ds $$

$$=    { (-1)^m\  \alpha^{-m}  \over 2\pi i}     \int_{\gamma - 2(k+1) -i\infty}^{\gamma- 2(k+1) +i\infty}\  (-s)_{m} \Gamma(s) \zeta (s) \zeta( s + 2k+1)\ (2\beta)^{-s}  ds $$

$$=   { (-1)^m\ \alpha^{-m}  \over 2\pi i}     \int_{\gamma - 2(k+1) -i\infty}^{\gamma- 2(k+1) +i\infty}\  {\Gamma(m-s)\over \Gamma (-s)} \Gamma(s) \zeta (s) \zeta( s + 2k+1)\ (2\beta)^{-s}  ds $$

$$=  {(-1)^{m+1}  \over 2\pi i}  \left({\beta\over \alpha}\right)^{m} {d^m\over d\beta^m}  \beta  {d \over d\beta} \int_{\gamma - 2(k+1) -i\infty}^{\gamma- 2(k+1) +i\infty}\  { \left[\ \Gamma(s)\ \right]^3  \over \Gamma (s+m) \Gamma (s+1-m)}  \zeta (s) \zeta( s + 2k+1)\ (2\beta)^{-s}  ds$$

$$=  {(-1)^{m+1}  \over 2\pi i}  \left({\beta\over \alpha}\right)^{m} \left(\beta {d\over d\beta} + m \right) {d^m\over d\beta^m}  \int_{\gamma - 2(k+1) -i\infty}^{\gamma- 2(k+1) +i\infty}\  { \left[\ \Gamma(s)\ \right]^3  \over \Gamma (s+m) \Gamma (s+1-m)}  \zeta (s) \zeta( s + 2k+1)\ (2\beta)^{-s}  ds.\eqno(6.16)$$
Now we employ the Mellin transform for the Legendre polynomials (cf. Entry 8.4.30.1 in [9], Vol. I) and the generalized Mellin-Parseval equality [14] to write from (6.16)

$$  {d^m\over d\alpha^m}  \sum_{n=1}^\infty {n^{-2k-1} \over  e^{2\beta n} -1 }  =  (-1)^{m+1}  \left({\beta\over \alpha}\right)^m \left(\beta {d\over d\beta} + m\right) {d^m\over d\beta^m} \int_{\beta}^\infty P_{m-1} \left( {2\beta\over t} -1 \right) \sum_{n=1}^\infty {n^{-2k-1} \over  e^{2 t n} -1 }  {dt\over t}.$$
Hence,  fulfilling consecutively  the differentiation under the integral sign due to the absolute and uniform convergence,  we find the chain of equalities 

$${d^m\over d\beta^m} \int_{\beta}^\infty P_{m-1} \left( {2\beta\over t} -1 \right) \sum_{n=1}^\infty {n^{-2k-1} \over  e^{2 t n} -1 }  {dt\over t} = {d^{m-1}\over d\beta^{m-1}} \bigg[   - {1 \over \beta} \sum_{n=1}^\infty {n^{-2k-1} \over  e^{2 \beta n} -1 }  + 2 \int_{\beta}^\infty P^\prime_{m-1} \left( {2\beta\over t} -1 \right) \sum_{n=1}^\infty {n^{-2k-1} \over  e^{2 t n} -1 }  {dt\over t^2}\bigg]$$

$$= - \sum_{q=0}^{m-1}  2^q P_{m-1}^{(q)} (1) \  {d^{m-1-q}\over d\beta^{m-1-q}} \bigg[ {1\over \beta^{q+1} } \sum_{n=1}^\infty {n^{-2k-1} \over  e^{2 \beta n} -1 } \bigg] $$

$$=  - \sum_{q=0}^{m-1} (-1)^{m-1-q} {2^q\over q!} \  P_{m-1}^{(q)} (1) \  \sum_{r=0}^{m-1-q}  (-1)^r (m-r-1)! \  \beta^{r-m} \ \binom{m-1-q}{r} {d^{r}\over d\beta^{r}}  \sum_{n=1}^\infty {n^{-2k-1} \over  e^{2 \beta n} -1 }$$

$$=  \sum_{q=0}^{m-1} (-1)^{q+1} {2^{m-1-q}\over (m-1-q)!} \  P_{m-1}^{(m-1-q)} (1) \  \sum_{r=0}^{q}  (-1)^r (m-r-1)! \  \beta^{r-m} \ \binom{q}{r} {d^{r}\over d\beta^{r}}  \sum_{n=1}^\infty {n^{-2k-1} \over  e^{2 \beta n} -1 }$$

$$=  - \sum_{r=0}^{m-1}  (m-r-1)! \  \beta^{r-m} \  {d^{r}\over d\beta^{r}} \left[  \sum_{n=1}^\infty {n^{-2k-1} \over  e^{2 \beta n} -1 } \right]  \sum_{q=0}^{m-1-r} (-1)^{q} {2^{m-1-r-q}\over (m-1-r-q)!} \  P_{m-1}^{(m-1-r-q)} (1) \  \binom{q+r}{r} $$

$$=  - \sum_{r=0}^{m-1}  r! \  \beta^{-1-r} \  {d^{m-1-r}\over d\beta^{m-1-r}} \left[  \sum_{n=1}^\infty {n^{-2k-1} \over  e^{2 \beta n} -1 } \right]  \sum_{q=0}^{r} (-1)^{q} {2^{r-q}\over (r-q)!} \  P_{m-1}^{(r-q)} (1) \  \binom{q+m-1-r}{q} .$$
But since (cf. [11])

$$  P_{m-1}^{(r-q)} (1) = {(r+m-1-q)! \over 2^{r-q} (r-q)! (m-1-r+q)!},$$
we get

$${d^m\over d\beta^m} \int_{\beta}^\infty P_{m-1} \left( {2\beta\over t} -1 \right) \sum_{n=1}^\infty {n^{-2k-1} \over  e^{2 t n} -1 }  {dt\over t} $$

$$=  (m-1)! \sum_{r=0}^{m-1}  {(-1)^{r+1}  \beta^{-1-r}\over  (m-1-r)!} \  {d^{m-1-r}\over d\beta^{m-1-r}} \left[  \sum_{n=1}^\infty {n^{-2k-1} \over  e^{2 \beta n} -1 } \right]  \sum_{q=0}^{r} (-1)^{q} \binom{r}{q} \ \binom{q+m-1}{m-1}.\eqno(6.17)$$
The inner sum is calculated in [9], Vol. I,  Entry 4.2.5.55, and it yields

$$\sum_{q=0}^{r} (-1)^{q} \binom{r}{q} \ \binom{q+m-1}{m-1} = (-1)^r \binom{m-1}{r}.$$
Therefore

$${d^m\over d\beta^m} \int_{\beta}^\infty P_{m-1} \left( {2\beta\over t} -1 \right) \sum_{n=1}^\infty {n^{-2k-1} \over  e^{2 t n} -1 }  {dt\over t} $$

$$=  - \sum_{r=0}^{m-1}  \beta^{r-m} \  (m-1-r)! \  \binom{m-1}{r}^2  \  {d^{r}\over d\beta^{r}} \left[  \sum_{n=1}^\infty {n^{-2k-1} \over  e^{2 \beta n} -1 } \right],$$ 
and, accordingly, (see (6.17)),

$$\alpha^m  {d^m\over d\alpha^m}  \sum_{n=1}^\infty {n^{-2k-1} \over  e^{2\beta n} -1 }  =    (- 1)^m \sum_{r=0}^{m-1}  \beta^{r} \  (m-1-r)! \  \binom{m-1}{r}^2  \bigg[ r+ \beta {d\over d\beta} \bigg]   \  {d^{r}\over d\beta^{r}} \left[  \sum_{n=1}^\infty {n^{-2k-1} \over  e^{2 \beta n} -1 } \right].$$ 
Recalling the Ramanujan formula (1.29), we write the left-hand side of the previous equality as follows

$$ \alpha^m  {d^m\over d\alpha^m} \sum_{n=1}^\infty {n^{-2k-1} \over  e^{2\beta n} -1 }  = \alpha^m  {d^m\over d\alpha^m} \bigg[  (-1)^k \pi^{2k} \alpha^{-2k} \ {\zeta(2k+1)\over 2}  +   2^{2k} \pi^{2k+1} \sum_{q=0}^{k+1} (-1)^{k+q}  {B_{2q} B_{2(k+1-q)} \pi^{2q-1} \over  (2q)! (2(k+1-q))!} \  \alpha^{1-2q} \bigg.$$

$$\bigg.+  (-1)^k \pi^{2k} \alpha^{-2k}  \sum_{n=1}^\infty {n^{-2k-1} \over  e^{2\alpha n} -1 }\bigg]$$

$$=  (-1)^{k+m}  {(2k-1+m)!\over 2 (2k-1)!} \left({\pi\over \alpha}\right)^{2k} \  \zeta(2k+1) +    2^{2k} \pi^{2k+1} \sum_{q=1}^{k+1} (-1)^{k+q+m}  {B_{2q} B_{2(k+1-q)}  (2(q-1)+m)! \over  (2q)! (2(q-1))! (2(k+1-q))!} \  \left({\pi\over \alpha}\right)^{2q-1} \bigg.$$

$$\bigg. + {(-1)^{k}\over (1-m)!}  2^{2k} \pi^{2k}\alpha {B_{2(k+1)} \over  (2(k+1))!}  +  (-1)^k \pi^{2k} \alpha^m  {d^m\over d\alpha^m} \bigg[ \alpha^{-2k}  \sum_{n=1}^\infty {n^{-2k-1} \over  e^{2\alpha n} -1 }\bigg]$$

$$=  (-1)^{k+m}  {(2k-1+m)!\over 2 (2k-1)!} \left({\pi\over \alpha}\right)^{2k} \  \zeta(2k+1) +    2^{2k} \pi^{2k+1} \sum_{q=1}^{k+1} (-1)^{k+q+m}  {B_{2q} B_{2(k+1-q)}  (2(q-1)+m)! \over  (2q)! (2(q-1))! (2(k+1-q))!} \  \left({\pi\over \alpha}\right)^{2q-1} \bigg.$$

$$\bigg.+ {(-1)^{k}\over (1-m)!}  2^{2k} \pi^{2k}\alpha {B_{2(k+1)} \over  (2(k+1))!} + {(-1)^k \over  (2k-1)!}  \left({\pi\over \alpha}\right)^{2k}    \sum_{q=0}^{m} (-1)^{m+q} \  (2k-1+m-q)!  \binom{m}{q} \alpha^{q}  {d^q\over d\alpha^q}  \sum_{n=1}^\infty {n^{-2k-1} \over  e^{2\alpha n} -1 }.$$
Thus we establish the equality 

$${\beta^{m}\over m!}   {d^{m}\over d\beta^{m}}  \sum_{n=1}^\infty {n^{-2k-1} \over  e^{2 \beta n} -1 } +    \sum_{r=1}^{m-1}  {r\over (r-1)!} \ \binom{m-1}{r}  \ \beta^{r} \  {d^{r}\over d\beta^{r}}  \sum_{n=1}^\infty {n^{-2k-1} \over  e^{2 \beta n} -1 } $$

$$=  (-1)^{k}  \binom{2k-1+m}{m}  \left({\pi\over \alpha}\right)^{2k} \bigg[ {\zeta(2k+1)\over 2} + \sum_{n=1}^\infty {n^{-2k-1} \over  e^{2\alpha n} -1 } \bigg] + {(-1)^{k+m}\over (1-m)!}  2^{2k} \pi^{2k}\alpha {B_{2(k+1)} \over  (2(k+1))!} $$

$$ +    2^{2k} \pi^{2k+1} \sum_{q=1}^{k+1} (-1)^{k+q}  {B_{2q} B_{2(k+1-q)} \over  (2q)!  (2(k+1-q))!} \  \binom{2(q-1)+m}{m}  \left({\pi\over \alpha}\right)^{2q-1} $$

$$+ {(-1)^{k+m} \over m!}  \left({\pi\over \alpha}\right)^{2k} \ \alpha^{m}  {d^m\over d\alpha^m}  \sum_{n=1}^\infty {n^{-2k-1} \over  e^{2\alpha n} -1 } +  (-1)^k \  \left({\pi\over \alpha}\right)^{2k}    \sum_{r=1}^{m-1} {(-1)^{r}\over r!}  \binom{2k-1+m-r}{2k-1} $$

$$\times  \alpha^{r}  {d^r\over d\alpha^r}  \sum_{n=1}^\infty {n^{-2k-1} \over  e^{2\alpha n} -1 },\quad m \in \mathbb{N}$$
which is equivalent to the Ramanujan-type formula  (1.30) and its particular cases (1.33), (1.34), (1.35), (1.36).  We invite interested readers to derive possibly new particular cases, comparing Ramanujan-type formulae (1.30), (6.8), (6.12), (6.14).

Finally, let us  extend the Mellin transform method to  the Mellin-Barnes type integrals, involving the Dirichlet beta function.   Indeed, let $a, b > 0, ab= \pi^2/4$. Using the reciprocal Mellin-Barnes representation for Entry 2.4.3.2 in [9], Vol. I and interchanging the order of integration and summation due to the dominated convergence, we deduce the equality

$$ \sum_{n=0}^\infty {(-1)^n \   (2n+1)^{- 2k-1}\over  \cosh((2n+1)a)} = {1\over \pi i}  \int_{\mu -i\infty}^{\mu +i\infty}\   \Gamma(s) \beta (s) \beta(s+2k+1) a^{-s} ds, \quad \mu > 1.\eqno(6.18)$$
Then functional equation (1.4) for the Dirichlet beta function yields

$$  \sum_{n=0}^\infty {(-1)^n \   (2n+1)^{- 2k-1}\over  \cosh((2n+1)a)} = {1\over \pi i}  \int_{\mu -i\infty}^{\mu +i\infty}\   \Gamma(s) \beta (s) \beta(s+2k+1) a^{-s} ds$$

$$= {(-1)^{k+1} \over \pi i} \ \left({\pi\over 2}\right)^{2k-1}  \int_{\mu -i\infty}^{\mu +i\infty}\   \Gamma(s) \left({\pi\over 2}\right)^{2s} \cos(\pi s/2) \sin(\pi s/2) \Gamma(1-s) \beta (1-s) $$

$$\times \Gamma(-s-2k) \beta(-s-2k) a^{-s} ds = {(-1)^{k+1} \pi \over 2 \pi i} \ \left({\pi\over 2}\right)^{2k-1}  \int_{\mu -i\infty}^{\mu +i\infty}\   \beta (1-s) \Gamma(-s-2k) \beta(-s-2k)  b^{s} ds $$

$$= {(-1)^{k+1}  \over \pi i}\  \left({a\over  b }\right)^{k}  \int_{-2k-\mu -i\infty}^{-2k-\mu +i\infty}\    \Gamma(s) \beta(s) \beta (s+2k+1) b^{-s} ds.$$
Assuming $1 < \mu < 2$, we shift the contour to the right, taking into account residues at simple left-hand poles of the gamma function $s= - 2 m,\  m=0 ,1,\dots, k.$
Then by residue theorem we obtain the formula 

$$     \sum_{n=0}^\infty {(-1)^n \   (2n+1)^{- 2k-1}\over  \cosh((2n+1)a)}  =    \left({a \over b}\right)^{k} \sum_{n=0}^\infty {(-1)^{k+n+1} \   (2n+1)^{- 2k-1}\over  \cosh((2n+1)b)}$$

$$ + {1 \over 2}  \left({\pi \over 2}\right)^{2k+1} \sum_{m=0}^k { (-1)^m E_{2m} E_{2(k-m)} \over (2m)! (2(k-m)!} \ \left({b\over a}\right)^{m-k},\quad k \in \mathbb{N}_0,\quad  ab= {\pi^2\over 4}.\eqno(6.19)$$
Particular cases of the latter identity are confirmed by Entries 5.3.7.1, 5.3.7.3 in [9], Vol. I.   

Further,  using   the reciprocal Mellin-Barnes representation for  Entry 2.3.14.2 in [9], Vol. I, we have  

$${1\over e^{cx}+1} = {1\over 2 \pi i}     \int_{\mu -i\infty}^{\mu +i\infty}\  \left( 1- 2^{1-s}\right) \Gamma(s) \zeta (s)  \left( cx\right)^{-s} ds,\ \mu > 1, \ c >0,$$
and then interchange  the order of summation and integration to  derive 

$$   \sum_{n=0}^\infty {(-1)^n \ e^{2 (2n+1) a } \over (2n+1)^{2k} \left[ \exp(2(2n+1) a ) + 1\right]}$$

$$= \beta (2k) -  {1\over 2\pi i}     \int_{\mu -i\infty}^{\mu +i\infty}\  \left( 1- 2^{1-s}\right) \Gamma(s) \zeta (s)  \beta( s + 2k)  (2a)^{-s} ds.$$
Let $a,b > 0, \ ab= \pi^2$. Then  functional equations for the Dirichlet beta function (1.4) and Riemann zeta function (1.6) yield the equalities 

$$   \sum_{n=0}^\infty {(-1)^n \ e^{2 (2n+1) a } \over (2n+1)^{2k} \left[ \exp(2(2n+1) a ) + 1\right]}$$

$$   =  \beta (2k) -  {2^{1-2k}  \pi^{2k} \over 2\pi i}     \int_{\mu -i\infty}^{\mu +i\infty}\  \left( 1- 2^{1-s}\right) \Gamma(s) 
 \pi^{2(s-1)} \sin(\pi s/2) \cos(\pi (s+2k)/2) $$
 
 $$\times \Gamma(1-s) \Gamma(1-s-2k)\zeta (1-s)  \beta(1- s - 2k)  (2a)^{-s} ds$$

$$= \beta (2k) -  {2^{-2k}  \pi^{2k-1} (-1)^k \over 2\pi i}     \int_{\mu -i\infty}^{\mu +i\infty}\  \left( 1- 2^{1-s}\right)   \Gamma(1-s-2k)\zeta (1-s)  \beta(1- s - 2k)  2^{-s} b^s ds$$

$$= \beta (2k) +  (-1)^k  \left({a\over b}\right)^{k-1/2}  {1\over 4\pi i}   \int_{1-2k-\mu -i\infty}^{1-2k-\mu +i\infty}\  \left(  1-2^{-s-2k}\right)   \Gamma(s)\beta(s) \zeta (s+2k) 2^s b^{-s}  ds.\eqno(6.20)$$
Thus  via the residue theorem we obtain 

$$  \sum_{n=0}^\infty {(-1)^{n+1} \  \over (2n+1)^{2k}  \left[ \exp(2(2n+1) a ) + 1\right]} =   (-1)^k  2^{-2k}  \left({a\over b}\right)^{k-1/2}  {1\over 4\pi i}   \int_{1-2k-\mu -i\infty}^{1-2k-\mu +i\infty}\  \left(  2^{s+2k}-1 \right)    \Gamma(s)\beta(s) \zeta (s+2k)  b^{-s}  ds$$

$$=  {(-1)^k \over 2}  \left({a\over b}\right)^{k-1/2}  \sum_{n=0}^\infty {(-1)^{n} \  (2n+1)^{-2k} \over  \cosh((n+1/2)b)}- {1\over 2} \  \beta(2k) + (- 1)^{k}  4^{-k-1}   \pi^{2k}  \sum_{m=1}^{k}  (-1)^{m}   \left(  4^{m} -1\right)  {B_{2m}  E_{2(k-m)}  4^m   \over (2m)!  (2(k-m))!}  \left({a\over b}\right)^{m-1/2},$$
proving a companion identity for (6.19)

$$  \sum_{n=0}^\infty  (-1)^{n} \  { \tanh((2n+1) a)  \over (2n+1)^{2k} }  =  (-1)^k \left({a\over b}\right)^{k-1/2}  \sum_{n=0}^\infty {(-1)^{n} \  (2n+1)^{-2k} \over  \cosh((n+1/2)b)} $$

$$+ {(- 1)^{k}\over 2}  \left({\pi\over 2}\right)^{2k} \sum_{m=1}^{k}  (-1)^{m}   \left(  4^{m} -1\right) \ 4^m \  {B_{2m} \ E_{2(k-m)}   \over (2m)!  (2(k-m))!}  \left({a\over b}\right)^{m-1/2},\   k \in \mathbb{N},  \ a,b >0,\   ab= \pi^2.$$
In the meantime, Entry 2.3.14.6 in [9], Vol. I suggests the Mellin-Barnes representation for the series 

$$   \sum_{n=0}^\infty {(-1)^n (2n+1)^{-2k} \over  \exp((2n+1) a ) - 1 } =  {1\over 2\pi i}     \int_{\mu -i\infty}^{\mu +i\infty}\  \Gamma(s) \zeta (s)  \beta( s + 2k)  a^{-s} ds.$$
Then in the same manner as in (6.20) we get

$$   \sum_{n=0}^\infty {(-1)^n (2n+1)^{-2k} \over  \exp((2n+1) a ) - 1 } =    (-1)^k 2^{-2k}  \left({a\over b}\right)^{k-1/2}  {1\over 2\pi i}    \int_{1-2k-\mu -i\infty}^{1-2k-\mu +i\infty}\    \Gamma(s)  \beta( s ) \zeta (s+2k)   b^{-s} ds$$

$$ = - {1\over 2}\ \beta (2k) +  (-1)^k \  2^{-2k-1} \left({ a\over b}\right)^{k-1/2} \sum_{n=1}^\infty  {n^{-2k} \over \cosh( b n) }   $$

$$ +  {a^{2k}\over 4}  \sum_{m=0}^k { (-1)^{m}  E_{2m} B_{2(k-m)} \over  4^m (2m)! (2(k-m))!} \  \left({b\over a}\right)^{m+1/2}$$
which yields the following   Ramanujan  formula (see [1], Ref. [10], p. 277, Entry 21(iii))

$$   \sum_{n=0}^\infty {(-1)^n (2n+1)^{-2k} \over  \exp((2n+1) a ) - 1 } = - {1\over 2}\ \beta (2k) +  (-1)^k \  2^{-2k-1} \left({ a\over b}\right)^{k-1/2} \sum_{n=1}^\infty  {n^{-2k} \over \cosh( b n) }   $$

$$ +  {a^{2k}\over 4}  \sum_{m=0}^k { (-1)^{m}  E_{2m} B_{2(k-m)} \over  4^m (2m)! (2(k-m))!} \  \left({b\over a}\right)^{m+1/2} ,\   k \in \mathbb{N},  \ a,b >0,\   ab= \pi^2.$$

\bigskip
\centerline{{\bf Acknowledgments}}
\bigskip

\noindent The work was partially supported by CMUP, which is financed by national funds through FCT (Portugal)  under the project with reference UIDB/00144/2020.

\bigskip
\centerline{{\bf References}}
\bigskip
\baselineskip=12pt
\medskip
\begin{enumerate}

\item[{\bf 1.} \ ]  B.C. Berndt,  A. Straub,  Ramanujan's formula for $\zeta(2n+1)$. {\it Exploring the Riemann zeta function, 13Ð34, Springer, Cham}, 2017.

\item[{\bf 2.} \ ]  E. Grosswald, Die Werte der Riemannschen Zeta-funktion an ungeraden Argumentstellen, {\it Nachr. Akad. Wiss. G{\" o}ttingen} (1970), 9-13.

\item[{\bf 3.} \ ]  A.P. Guinand, Functional equations and self-reciprocal functions connected with Lambert series, {\it Quart. J. Math.} {\bf 15}  (1944), 11-23.

\item[{\bf 4.} \ ]    A.P. Guinand, Some rapidly convergent series for the Riemann $\xi$-function, {\it Quart. J. Math. ser. (2)},  {\bf 6} (1955), 156-160.

\item[{\bf 5.} \ ]   S.Kanemitsu,Y.Tanigawa, and M.Yoshimoto, On rapidly convergent series for the Riemann zeta-values via the modular relation,  {\it Abh. Math. Sem. Univ. Hamburg } {\bf 72} (2002), 
187- 206.

\item[{\bf 6.} \ ]  K.S. K{\" o}lbig,  The polygamma function $\psi^{(k)} (x)$ for $x= 1/4$ and $x= 3/4$. {\it Journal of Computational and Applied Mathematics }, {\bf 75} (1996), 43-46.

\item[{\bf 7.} \ ]  A.F. Loureiro,  P. Maroni and S. Yakubovich,  On a polynomial sequence associated with the Bessel operator. {\it Proceedings of the AMS}, {\bf 142} (2014), N2,  467-482.

\item[{\bf 8.} \ ]  S.L. Malurkar, On the application of Herr MellinÕs integrals to some series, {\it J. Indian Math. Soc.} {\bf 16}  (1925/26), 130Ð138.

\item[{\bf 9.}\ ] A.P. Prudnikov, Yu.A. Brychkov and O.I. Marichev, {\it Integrals and Series}. Vol. I: {\it Elementary
Functions}, Vol. II: {\it Special Functions}, Gordon and Breach, New York and London, 1986, Vol. III : {\it More special functions},  Gordon and Breach, New York and London,  1990.

\item[{\bf 10.}\ ]  J. Riordan, Combinatorial Identities, Wiley, New York, 1968.

\item[{\bf 11.}\ ]  D.E. Winch, P.H. Roberts, Derivatives of addition theorems for Legendre functions. {\it J. Austral. Math. Soc. Ser. B.}, {\bf 37} (1995),  212-234.

\item[{\bf 12.}\ ]  S. Yakubovich,  Certain identities, connection and explicit formulas for the Bernoulli and Euler numbers and the Riemann zeta-values. {\it Analysis (Berlin)}, {\bf 35} (2015), N 1, 59-71. 

\item[{\bf 13.}\ ]  S. Yakubovich, {\it  Index Transforms}, World Scientific Publishing Company, Singapore, New Jersey, London and Hong Kong, 1996.

\item[{\bf 14.}\ ] S.  Yakubovich and Yu.  Luchko, {\it The Hypergeometric Approach to Integral Transforms and Convolutions},  Mathematics and Applications, Vol. 287, Kluwer Academic Publishers, Dordrecht, 1994.

\end{enumerate}

\end{document}